\documentclass[10pt]{article}
\usepackage{amssymb} 
\newcommand{\be}{\begin{equation}}
\newcommand{\ee}{\end{equation}}
\newcommand{\bea}{\begin{eqnarray}}
\newcommand{\eea}{\end{eqnarray}}
\newcommand{\nn}{\nonumber \\}
\newcommand{\cD}{\mathbb{D}}
\newcommand{\cZ}{\mathbb{Z}}
\newcommand{\cP}{\mathbb{P}}
\newcommand{\cN}{\mathbb{N}}
\newcommand{\Col}{{\,\rm Col}}
\newcommand{\cR}{\mathbb{R}}
\newcommand{\Sir}{{\,\rm Syr}}
\newcommand{\Syr}{{\bf Syrac}}
\newcommand{\Osc}{{\,\rm Osc}}
\newcommand{\nicht}[1]{ }
\begin{document}

\title{The Collatz Problem generalized to $3x+k$}

\author{Franz Wegner \\
Institut f\"ur Theoretische Physik \\
Universit\"at Heidelberg \\
D-69120 Heidelberg Germany}

\maketitle

\begin{abstract}
The Collatz problem with $3x+k$ is revisited. Positive and
negative limit cycles are given up to $k=9997$ starting with
$x_0=-2\cdot10^7...+2\cdot10^7$.
A simple relation between the probability distribution for the
Syracuse iterates for various $k$ (not divisible by 2 and 3)
is obtained. From this it follows, that the oscillation
considered by Tao 2019 (arXiv:1909.03562v2) does not depend
on $k$. Thus this piece of the proof of his theorem 1.3 
"Almost all Collatz orbits attain almost bounded values" 
holds for all $k$ not divisible by 2 and 3.
\end{abstract}

\section{Introduction}

\paragraph{Collatz problem and Collatz conjecture}
The Collatz problem starts with a positive integer $x_0$
and defines a sequence $x_1,x_2,...$ by the rule
\be
x_{i+1} = \Col(x_i), \quad
\Col(x) := \left\{ \begin{array}{cc} 3x+1, & x \mbox{ odd}, \\
x/2 & x \mbox{ even} \end{array} \right.
\ee
The Collatz conjecture is: For all positive $x_0$ the sequence
tends to the limit cycle $C_1$.
\be
C_1=1,4,2,1,... \label{Cp1}
\ee
Starting with negative $x_0$ three limit cycles are known:
\bea
C_{-1} &=& -1,-2,-1,..., \\
C_{-5} &=& -5,-14,-7,-20,-10,-5,..., \\
C_{-17}&=& -17,-50,-25,-74,-37,-110,-55,-164,-82,-41,-122, \nn
&& -61,-182,-91,-272,-136,-68,-34,-17,... \label{Cm17}
\eea
\paragraph{Generalized Collatz problem}
Here we consider a generalization of the Collatz problem.
We replace $3n+1$ by $3n+k$ with odd $k$. Thus this generalized
map is defined by the rule
\be
x_{i+1} = \Col_k(x_i), \quad
\Col_k(x) := \left\{ \begin{array}{cc} 3x+k, & x \mbox{ odd}, \\
x/2 & x \mbox{ even} \end{array} \right.
\ee
The Collatz-problem uses $k=1$.

The question is, what happens to the sequence $x_i$ for large $i$.
There are two possibilities: The sequence may end in a limit cycle,
or the sequence may tend to $\pm\infty$.

The Collatz problem has attracted many mathematicians.
Lagarias has collected many papers\cite{Lag03,Lag06}
on this problem and published them with short useful comments.
The $3x+k$-problem is sometimes hidden as $3x+1$-problem
by allowing $x$ to be rational with denominator $x\in\cD$
\cite{Bre02,Lag90} by use of the second equation (\ref{comfac}),
where $\cD$ is the set of integers $x$ with $x=\pm1\bmod6$.

In section 2 we review properties of limit cycles
and present numerical results for limit cycles for $k=1$ to
$k=149$ in the tables \ref{tb1} to \ref{tb2}.
The limit cycles are given from $k=1$ to $k=9997$ in an
accompanying supplement\cite{Weglist}.
These are the cycles obtained by starting from $x_0=-2\cdot10^7$
to $x_0=+2\cdot10^7$. (To be precise: For reasons explained
at the beginning of section 2 only $k,x_0\in\cD$ were
considered). Comparable numerical calculations were performed
by Belaga and Mignotte\cite{BelMig06a,BelMig06b,BelMig12}
for $k=1$ to $k=19997$, but only for positive $x$.
All sequences ended in limit cycles. Thus I conjecture that there is
no unbounded sequence for odd $k$. Moreover it seems likely that
for given $k$ the number of limit cycles is finite, but it may
be very large. These are conjectures already expressed by
Lagarias\cite{Lag90}.

In section 3 I compare two quantities considered by Tao in his
paper\cite{Tao} on the Collatz problem.
One quantity is the probability distribution of
the Syracuse iterates. A simple connection between
those for $k\in\cD$ and $k=1$ is given.
Due to this relation the oscillation $\Osc_{m,n}$
is equal for all $k\in\cD$. Thus the oscillation has the same
decay property for all these $k$. This is an essential
piece of the proof of his theorem 1.3: {\it Almost all
Collatz orbits attain almost bounded values}.

\section{Limit Cycles}

In the first part of this section some properties of Collatz and 
Syracuse mapping are reviewed. Compare
\cite{Lag90,Lag03,Lag06,BelMig06a,BelMig06b,BelMig12}.
Then explicit limit cycles are given.

As mentioned before the set of integers $x$ with $x=\pm 1\bmod 6$
is denoted by $\cD$.
Thus even integers and multiples of 3 are excluded. Later we will use
\paragraph{Lemma 1} Let $f,g,h\in\cZ$ and $fg=h$.
If $f,g\in\cD$, then $h\in\cD$. If $h\in\cD$, then $f,g\in\cD$.

\paragraph{Choice of $k$} It is sufficient to restrict to
$k\in\cD$ for the following reasons:\\
(i) If one choose an even $k$, then as soon, as $x_i$ is odd,
all $x_j$ with $j>i$ will also be odd, so that $x_j$ will run
to $\pm\infty$, unless $x_i=-a/2$, in which case $x_j$ stays
constant.\\
(ii) If $k$ is a multiple of 3, $k=3k'$, then after one step also
$x_{i+1}=3x'_{i+1}$ with $x'_{i+1}=x_i+k'$ is a multiple of 3
and all $x_j=3x'_j$ for $j>i$ with $x'_{j+1}=3x'_j+k'$.
Thus any limit cycles for $k$ are known from from those of $k'$.

\paragraph{Choice of $x_0$} It is also sufficient to consider
only $x_0\in\cD$ for the following reasons:\\
(i) If $x_0$ is even, then it will be divided by two until an
odd number appears in the sequence. Thus it is sufficient to
consider only odd $x_0$.\\
If $x_i$ is a multiple of 3 and odd, then $x_{i+1}$ will not be
a multiple of 3, provided that $k\in\cD$. From then on there
will be no more $x_j$ divisible by 3.

\paragraph{Common factors of $x$ and $k$}
If $x$ and $k$ have the common factor $p$,
\be
x=px', \quad k=pk', \quad x,k \in\cD, \quad p\in\cZ
\ee
then
\be
x',k',p\in\cD, \quad \Col_k(x) = p \Col_{k'}(x'). \label{comfac}
\ee
This has two consequences:\\
(i) If $x$ and $k$ have a common factor $p$, then it is sufficient to
follow $x'$ with $\Col_{k'}$.
If $x_0=px'_0$ and $k=pk'$ with integer $p$, then
the sequences generated from $x_i$ with $k$ and $x'_i$ with
$k'$ obey $x_j=px'_j$.\\
(ii) It is sufficient to consider only sequences for positive $k$, but
for positive and negative $x$, since one obtains the sequences
for negative $k$ from those of positive $k$ with $p=-1$.

\paragraph{Syracuse formulation} In this formulation the step
$3x+k$ and any following divisions by two are put into one
step denoted by $\Sir$,
\be
\Sir_k(x) = \frac{3x+k}{2^{\nu(3x+k)}}, \label{Sir1}
\ee
where $\nu(y)$ is the number of factors 2 in $y$, thus
$y/2^{\nu(y)}\in 2\cZ+1$.
Changing the indices $i$ of $x$, we consider now the Syracuse
mapping
\be
x_{i+1} = \Sir_k(x_i). \label{Sir2}
\ee
It maps the odd $x_i$ into the next odd $x$ of the Collatz map,
which is now  denoted by $x_{i+1}$.

Iterating eq. (\ref{Sir2}) one obtains
\be
2^{a_{[i,j-1]}} x_j = 3^{j-i} x_i + F_{j,i} k,
\ee
with the notation
\bea
a_{[i,j]} &:=& \sum_{l=i}^j a_l, \quad a_l:=\nu(x_l), \\
F_{j,i} &:=& \sum_{j=0}^{i-1} 3^{j-i-1-l} 2^{a_{[i,i+l-1]}} \nn
&=& 3^{j-i-1} + 3^{j-i-2}2^{a_i} + ... 
+2^{a_{[i,j-2]}}.
\eea
$F_{j,i}$ yields for the smallest $j-i$
\be
F_{i,i}=0, \quad F_{i+1,i}=1, \quad F_{i+2,i}=3+2^{a_i}.
\ee
One derives the recursion relations
\be
F_{j,i} = 2^{a_i} F_{j,i+1} +3^{j-i-1} = 3F_{j-1,i} +2^{a_{[i,j-2]}}.
\label{Frec}
\ee

\paragraph{Limit cycles}
A limit cycle obeys $x_{i+b}=x_i$, which implies $a_{i+b}=a_i$,
then
\be
D x_i= F_{i+b,i} k, \quad D := 2^a -3^b, \quad a=a_{[i,i+b-1]}.
\label{DF}
\ee
$a$ does not depend on $i$.
Thus with the greatest common divisor of $D$ and $F$
one obtains the solutions for limit cycles
\be
k = \frac{qD}{\gcd(D,F)}, \quad
x_i = \frac{qF_{i+b,i}}{\gcd(D,F)}.
\ee
Since $D,F\in\cD$, Lemma 1 yields $\gcd(D,F_{i+b,i})\in\cD$.
With $q\in\cD$ one obtains also $k,x_i\in\cD$. 
$F_{j,i}$ is positive for $j>i$. $D$ can be positive or negative.
We choose the sign of $\gcd(F,D)$ equal to the sign of $D$.
The recurrence relation (\ref{Frec}) yields for $j=i+b+1$
\be
2^{a_i} F_{i+b+1,i+1} = 3 F_{i+b,i} + D. \label{Frec2}
\ee
Call $g_i=\gcd(D,F_{i+b,i})$. Then the right hand side of
(\ref{Frec2}) is a multiple of $g_i$ and hence also the l.h.s.
One concludes that $F_{i+b+1,i+1}$ is a multiple of $g_i$,
since $g_i\in\cD$ does not contain factors 2.
Similarly since $F_{i+b+1,i+1}$ is a multiple of $g_{i+1}$,
also $3F_{i+b,i}$ and $F_{i+b,i}$ itself are multiples of
$g_{i+1}$. Consequently $g_{i+1}=g_i$ and thus
$\gcd(D,F_{i+b,i})$ does not depend on $i$, which allows
to skip the indices of $F$ in this greatest common divisor.

If the ordered list $\{a\}=(a_0,a_1,...,a_{b-1})$ contains subperiods,
\be
b=pb', \quad a_j:=a_{j+b'}, \quad a=pa', \quad p>1,
\ee
then also $F_{j+b',j}=F_{j,j-b'}$
and one obtains
\be
2^{a'} (x_{j+b'}-\gamma_j) = 3^{b'} (x_j-\gamma_j), \quad
\gamma_j=\frac{F_{j+b',j}}{D'}k, \quad D'=2^{a'}-3^{b'},
\label{xrecF}
\ee
Since $F_{j+b',j}$ and $\gamma_j$ are periodic in $j$ 
with period $b'$, iteration of (\ref{xrecF}) $p$ times
yields
\be
2^a (x_{j+b}-\gamma_j) = 3^b (x_j-\gamma_j).
\ee
Since $x_{j+b}=x_j$,one obtains $x_j=\gamma_j$ and thus
$x_{j+b'}=x_j$. Thus the shorter period $b'$ of $a_j$
implies also the period $b'$ of $x_j$.
Thus a periodic set $(a_0,a_1,...,a_{b-1})$ with subperiods
yields the same results as $(a_0,a_1,...,a_{b'-1})$.
Hence it is sufficient to consider lists $\{a\}$ without
periodic subsets. The corresponding cycles are called irreducible. 

The number of ordered partitions
$(a_0,a_1,...,a_{b-1})$ is
\be
\#(a,b)=\frac{(a-1)!}{(b-1)!/(a-b)!}.
\ee
The $b$ cyclic permutations of $a_i$ yield a limit cycle. 
If there are no periodic subcycles, which is the case,
if $a$ and $b$ are coprime, then the number of limit cycles is
\be
N_{cyc} = \frac{(a-1)!}{b!(a-b)!}. \label{Ncyc1}
\ee
If $a$ and $b$ are not coprime,
\be
\gcd(a,b) = \prod_i p_i^{\nu_i}, \quad \nu_i\ge 1,
\ee
then one has to subtract those with subperiods and obtains the
number of cycles without periodic subcycles
\be
N_{cyc} = \frac1b \sum (-)^{\sum_i n_i} 
\#(a/\prod_ip_i^{n_i},(b/\prod_ip_i^{n_i}). \label{Ncyc2}
\ee
The summation runs over all $n_i=0,1$ independently.

\paragraph{Table of limit cycles}
I determined the limit cycles for $k=1...9997,k\in\cD$
for starting $x_0=-2\cdot10^7...+2\cdot10^7$, $x_0\in\cD$,
$\gcd(x_0,k)=1$. The Collatz sequence ended in all these
cases in a limit cycle.
If $\gcd(k,x_0)=q\not=1$, then
$\Col_k^i(x_0)=q\Col_{q/k}^i(x_0/q)$
and thus the limit cycle of $x_0$ is already contained in the
limit cycles for $k'=k/q$. This holds in particular for q=k.

The limit cycles are listed in tables
\ref{tb1} to \ref{tb2} from $k=1$ up to $k=149$.
In separate accompanying tables\cite{Weglist} I give the limit cycles
for $k=1...9997:k\in\cD$. Belaga and Mignotte have calculated
them for $k=1...19999:k\in\cD$, but only for positive 
cycles.\cite{BelMig06b} Unfortunately I found only 9 introductory
pages of this paper in the internet. The table was missing.

Table 3.1 of Lagarias\cite{Lag90} gives lowest $x_{lc}$ and the
corresponding $a$ for $k=1$ to $149$. They agree with those
given in the tables \ref{tb1} to \ref{tb2}. $x_{lc}=7$ is not
found in the table, since $\gcd(x_{lc},k) = 7.$ Thus
the limit cycle is contained in that for $k'=5$.

\begin{table}

\noindent\begin{tabular}{rrrrr}
 $k$ & $x_0$ & $x_{lc}$ & $a$ & $b$ \\
\hline
  1 &     1 &     1 &    2 &    1 \\
  1 &    -1 &    -1 &    1 &    1 \\
  1 &    -5 &    -5 &    3 &    2 \\
  1 &   -17 &   -17 &   11 &    7 \\ \hline
  5 &     1 &     1 &    3 &    1 \\
  5 &     7 &    19 &    5 &    3 \\
  5 &    23 &    23 &    5 &    3 \\
  5 &   187 &   187 &   27 &   17 \\
  5 &   259 &   347 &   27 &   17 \\ \hline
  7 &     1 &     5 &    4 &    2 \\ \hline
 11 &     1 &     1 &    6 &    2 \\
 11 &     5 &    13 &   14 &    8 \\
 11 &   -19 &   -19 &    4 &    3 \\ \hline
 13 &     1 &     1 &    4 &    1 \\
 13 &    19 &   131 &   24 &   15 \\
 13 &   155 &   211 &    8 &    5 \\
 13 &   163 &   251 &    8 &    5 \\
 13 &   187 &   287 &    8 &    5 \\
 13 &   191 &   259 &    8 &    5 \\
 13 &   227 &   227 &    8 &    5 \\
 13 &   283 &   283 &    8 &    5 \\
 13 &   319 &   319 &    8 &    5 \\ \hline
 17 &     1 &     1 &    7 &    2 \\
 17 &    11 &    23 &   31 &   18 \\
 17 &   -41 &   -65 &    6 &    4 \\
 17 &   -73 &   -73 &    6 &    4 \\ \hline
 19 &     1 &     5 &   11 &    5 \\
 19 &   -67 &  -115 &   17 &   11 \\ \hline
 23 &     1 &    41 &   43 &   26 \\
 23 &     5 &     5 &    5 &    2 \\
 23 &     7 &     7 &    5 &    2 \\
 23 & -2263 & -2263 &   19 &   12 \\
 23 & -2359 & -2359 &   19 &   12 \\
 23 & -2495 & -2963 &   19 &   12 \\
 23 & -2503 & -3743 &   19 &   12 \\
 23 & -2567 & -3415 &   19 &   12 \\
 23 & -2743 & -2743 &   19 &   12 \\
 23 & -3091 & -3091 &   19 &   12 \\
 23 & -4819 & -4819 &   19 &   12 \\ \hline
\end{tabular}
\hspace{5mm}
\noindent\begin{tabular}{rrrrr}
 $k$ & $x_0$ & $x_{lc}$ & $a$ & $b$ \\
\hline
 25 &     1 &     7 &   16 &    8 \\
 25 &    17 &    17 &    8 &    4 \\
 25 &   -89 &  -113 &   34 &   22 \\ \hline
 29 &     1 &     1 &    5 &    1 \\
 29 &     5 &    11 &   17 &    9 \\
 29 &  2531 &  3811 &   65 &   41 \\
 29 &  7055 &  7055 &   65 &   41 \\
 29 &   -85 &  -109 &   12 &    8 \\ \hline
 31 &     1 &    13 &   23 &   12 \\
 31 &   -79 &   -95 &    9 &    6 \\ \hline
 35 &     1 &    13 &    8 &    4 \\
 35 &    11 &    17 &    8 &    4 \\ \hline
 37 &     1 &    19 &    6 &    3 \\
 37 &     7 &    29 &    6 &    3 \\
 37 &    23 &    23 &    6 &    3 \\ \hline
 41 &     1 &     1 &   20 &    8 \\ \hline
 43 &     1 &     1 &   11 &    3 \\ \hline
 47 &     1 &    25 &   28 &   16 \\
 47 &     5 &     5 &   18 &    7 \\
 47 &    37 &    65 &    7 &    4 \\
 47 &    41 &    85 &    7 &    4 \\
 47 &    53 &    89 &    7 &    4 \\
 47 &    73 &    73 &    7 &    4 \\
 47 &   101 &   101 &    7 &    4 \\ \hline
 49 &     1 &    25 &   38 &   22 \\
 49 &   -65 &   -65 &    5 &    4 \\ \hline
 53 &     1 &   103 &   29 &   17 \\ \hline
 55 &     1 &     1 &   12 &    4 \\
 55 &     7 &     7 &    6 &    2 \\
 55 &    41 &    41 &   28 &   16 \\ \hline
 59 &     1 &     1 &   28 &   11 \\
 59 &     5 &   133 &   10 &    6 \\
 59 &    61 &   181 &   10 &    6 \\
 59 &   125 &   217 &   10 &    6 \\
 59 &   149 &   149 &   10 &    6 \\
 59 &   169 &   185 &   10 &    6 \\
 59 &   221 &   221 &   10 &    6 \\ \hline
 61 &     1 &     1 &    6 &    1 \\
 61 &   175 &   235 &   66 &   41 \\ \hline
\end{tabular}
\caption{Limit cycles part 1}\label{tb1}
\end{table}

\begin{table}

\noindent\begin{tabular}{rrrrr}
 $k$ & $x_0$ & $x_{lc}$ & $a$ & $b$ \\
\hline
 65 &     1 &    19 &   24 &   12 \\ \hline
 67 &     1 &    17 &   30 &   16 \\ \hline
 71 &     1 &    29 &   10 &    5 \\
 71 &     7 &    31 &   10 &    5 \\
 71 &   893 &  4409 &   27 &   17 \\
 71 &  1073 &  2809 &   27 &   17 \\
 71 &  2585 &  2585 &   27 &   17 \\
 71 &  2633 &  3985 &   27 &   17 \\
 71 &  4121 &  4121 &   27 &   17 \\ \hline
 73 &     1 &    19 &   60 &   32 \\
 73 &     5 &     5 &   12 &    4 \\
 73 &     7 &    47 &   15 &    8 \\ \hline
 77 &     1 &     1 &   38 &   16 \\
 77 &  -397 & -1153 &   22 &   14 \\
 77 &  -685 &  -989 &   22 &   14 \\
 77 & -1165 & -1165 &   22 &   14 \\ \hline
 79 &     1 &     1 &   44 &   20 \\
 79 &     7 &     7 &   44 &   20 \\
 79 &   193 &   265 &   23 &   14 \\
 79 &   233 &   233 &   23 &   14 \\ \hline
 83 &     1 &   109 &   12 &    7 \\
 83 &    17 &   157 &   12 &    7 \\
 83 &    59 &    65 &   24 &   14 \\ \hline
 85 &     1 &     7 &  100 &   56 \\
 85 &  -293 &  -293 &   12 &    8 \\ \hline
 89 &     1 &    17 &   17 &    8 \\ \hline
 91 &     1 &     1 &   48 &   24 \\
 91 &    25 &    25 &   12 &    6 \\
 91 &    59 &    59 &   12 &    6 \\ \hline
 95 &     1 &     1 &   72 &   36 \\
 95 &     7 &    23 &   11 &    5 \\
 95 &    17 &    17 &   11 &    5 \\
 95 &  -511 &  -511 &   17 &   11 \\ \hline
 97 &     1 &     1 &   18 &    6 \\
 97 &     5 &    13 &    9 &    3 \\
 97 &  -449 &  -625 &   17 &   11 \\
 97 &  -505 &  -505 &   17 &   11 \\ \hline
101 &     1 &    11 &   14 &    6 \\
101 &     5 &    29 &    7 &    3 \\
101 &     7 &     7 &   14 &    6 \\
101 &    17 &    19 &    7 &    3 \\
101 &    23 &    23 &    7 &    3 \\
101 &    31 &    31 &    7 &    3 \\
\end{tabular}
\hspace{5mm}
\noindent\begin{tabular}{rrrrr}
 $k$ & $x_0$ & $x_{lc}$ & $a$ & $b$ \\
\hline
101 &    37 &    37 &    7 &    3 \\ \hline
103 &     1 &    23 &   45 &   22 \\
103 &     5 &     5 &   21 &    8 \\ \hline
107 &     1 &     1 &  106 &   53 \\ \hline
109 &     1 &    19 &   52 &   30 \\ \hline
113 &     1 &     1 &   29 &   12 \\
113 &  -769 &  -769 &   69 &   44 \\ \hline
115 &     1 &    13 &   27 &   13 \\
115 &     7 &    17 &   34 &   18 \\
115 &  -163 &  -211 &    7 &    5 \\
115 &  -227 &  -227 &    7 &    5 \\
115 &  -251 &  -251 &    7 &    5 \\ \hline
119 &     1 &     1 &   75 &   42 \\
119 &     5 &     5 &    7 &    2 \\
119 &    11 &    11 &    7 &    2 \\
119 &    13 &    23 &   39 &   18 \\
119 &   125 &   125 &   25 &   14 \\
119 &  -247 &  -311 &   18 &   12 \\ \hline
121 &     1 &     5 &   44 &   18 \\
121 &    19 &    19 &   88 &   46 \\
121 & -1085 & -1085 &   44 &   28 \\
121 & -1105 & -2089 &   22 &   14 \\ \hline
125 &     1 &     1 &    7 &    1 \\
125 &     7 &    47 &   33 &   19 \\
125 &   143 &   143 &   33 &   19 \\
125 &   431 &   899 &  118 &   74 \\ \hline
127 &     1 &     1 &   37 &   18 \\
127 &     5 &    41 &   37 &   18 \\ \hline
131 &     1 &    13 &   26 &   13 \\
131 &     7 &    17 &   26 &   13 \\
131 &    23 &    23 &   52 &   26 \\
131 & -7471 & -7471 &   38 &   24 \\
131 & -8003 & -15875 &   38 &   24 \\
131 & -10883 & -10883 &   38 &   24 \\
131 & -12107 & -16099 &   38 &   24 \\
131 & -13963 & -13963 &   38 &   24 \\ \hline
133 &     1 &    11 &   36 &   18 \\
133 &    31 &    59 &   42 &   24 \\ \hline
137 &     1 &     1 &   14 &    4 \\
137 &     5 &    41 &   72 &   40 \\
137 &    95 &   503 &   29 &   18 \\
137 &   415 &   967 &   29 &   18 \\
137 &   743 &   743 &   29 &   18 \\ \hline
\end{tabular}
\caption{Limit cycles part 2}\label{tb1a}
\end{table}

\begin{table}[t]

\noindent\begin{tabular}{rrrrr}
 $k$ & $x_0$ & $x_{lc}$ & $a$ & $b$ \\
\hline
139 &     1 &    11 &  136 &   74 \\
139 &  -331 & -1291 &   22 &   14 \\
139 &  -587 & -1931 &   22 &   14 \\
139 &  -779 & -2059 &   11 &    7 \\
139 &  -827 & -2203 &   22 &   14 \\
139 & -1043 & -1495 &   33 &   21 \\
139 & -1195 & -2251 &   11 &    7 \\
139 & -1339 & -2539 &   11 &    7 \\
139 & -1355 & -1355 &   22 &   14 \\
139 & -1367 & -1367 &   33 &   21 \\
139 & -1555 & -2971 &   11 &    7 \\
139 & -1835 & -2683 &   11 &    7 \\
139 & -1915 & -2507 &   11 &    7 \\
139 & -1979 & -2899 &   11 &    7 \\
139 & -2123 & -2123 &   11 &    7 \\
139 & -2131 & -2795 &   11 &    7 \\
139 & -2171 & -3187 &   11 &    7 \\
139 & -2195 & -3223 &   11 &    7 \\
139 & -2219 & -2219 &   11 &    7 \\
139 & -2347 & -2347 &   11 &    7 \\
139 & -2387 & -3511 &   11 &    7 \\
139 & -2419 & -2903 &   11 &    7 \\
139 & -2455 & -3227 &   11 &    7 \\
\end{tabular}
\hspace{5mm}
\begin{tabular}{rrrrr}
 $k$ & $x_0$ & $x_{lc}$ & $a$ & $b$ \\
\hline
139 & -2491 & -2491 &   11 &    7 \\
139 & -2519 & -2747 &   11 &    7 \\
139 & -2579 & -2579 &   11 &    7 \\
139 & -2603 & -2603 &   11 &    7 \\
139 & -2707 & -2707 &   11 &    7 \\
139 & -2711 & -2963 &   11 &    7 \\
139 & -2743 & -3287 &   11 &    7 \\
139 & -2939 & -2939 &   11 &    7 \\
139 & -3031 & -3031 &   11 &    7 \\
139 & -3155 & -3155 &   11 &    7 \\
139 & -3443 & -3443 &   11 &    7 \\
139 & -3479 & -3479 &   11 &    7 \\
139 & -3767 & -3767 &   11 &    7 \\ \hline
143 &     1 &     7 &  140 &   80 \\
143 &     5 &    29 &   16 &    7 \\
143 &    17 &    17 &   16 &    7 \\ \hline
145 &     1 &     1 &   10 &    2 \\
145 &     7 &    23 &   22 &   10 \\
145 &    17 &    47 &   34 &   18 \\
145 &  -209 &  -241 &   48 &   32 \\
145 &  -617 &  -617 &   12 &    8 \\ \hline
149 &     1 &    19 &   59 &   33 \\
149 &   395 &   667 &   42 &   26 \\ \hline
\end{tabular}

\caption{Limit cycles part 3}\label{tb2}
\end{table}

I give $k$, $a$, and $b$ for the limit cycles.
There are limit cycles with positive $x$ and those with
negative $x$. I call them positive and negative limit cycles.
For the positive limit cycles I list the smallest $x$ in
the cycle and denote it by $x_{lc}$. For negative limit cycles
$x_{lc}$ is that with the smallest $-x$ in the cycle.
Also $x_0$ is given. For positive
limit cycles it is the smallest positive $x_0$ from which the
cycle is reached, for negative cycles it is the smallest
$-x_0$ from which it is reached.

There is at least one positive limit cycle for each $k$.
For odd $x_i\in(-k/3,-1)$ the following $x_{i+1}$ is positive.
Thus it may happen that there is no negative limit cycle
under the condition $\gcd(k,x_0)=1$.
This is the case for $k=5,7,13$ and many more.
\begin{table}
\begin{tabular}{rrr}
$k$ & $N_+$ & $N_-$ \\ \hline
311 & 2 & 126 \\
3299 & 1 & 171 \\
5137 & 4 & 135 \\
6005 & 142 & 0 \\
6487 & 534 & 0 \\
6553 & 111 & 0 \\
7153 & 4 & 2908 \\
7463 & 162 & 0 \\
\end{tabular}
\hspace{5mm}
\begin{tabular}{rrr}
$k$ & $N_+$ & $N_-$ \\ \hline
7645 & 2 & 0 \\
7727 & 198 & 0 \\
7949 & 117 & 0 \\
8059 & 1 & 332 \\
8425 & 120 & 0 \\
8765 & 106 & 0 \\
9215 & 57 & 114 \\
9823 & 241 & 0
\end{tabular}
\caption{Number $N_+$ of positive and
$N_-$ of negative limit cycles for $k\in\cD$ up to 9997
with at least 100 limit cycles} \label{tb5}
\end{table}
\paragraph{Many limit cycles}
Table \ref{tb5} shows the number $N_+$ of positive and
$N_-$ of negative limit cycles for $k\in\cD$ up to 9997 with at least
100 limit cycles. In many cases $N_+$ is large and $N_-=0$
or $N_-$ is large and $N_+$ is small, typically $1...4$.
An exception is $k=9215$. Particularly many (2908) negative limit cycles
are obtained for $k=7153$.

\paragraph{Example $k=1$, $k=139$}
For $k=1$ there are the four limit cycles mentioned in eqs.
(\ref{Cp1} to \ref{Cm17}). The corresponding values of $D$ are
\bea
D_1=4-3=1, && D_{-1}=2-3=-1, \nn
D_{-5}=8-9=-1, && D_{-17}=2048-2187=-139.
\eea
The last cycle gives $D$ different from $|D|=1$. Inspection shows that
the corresponding $F$ are multiples of 139. Thus the cycle
contributes to $k=1$.
In the present case $(a,b)=(11,7)$ one has $N_{cyc}=30$. One appears
at $k=1$, the other 29 at  $k=139$. For $k=139$ only one positive
cycle with $(a,b)=(136,74)$ has been found. Besides the 29 negative
cycles with $(a,b)=(11,7)$ four negative cycles with $(a,b)=(22,14)$ and
two negative cycles with $(a,b)=(33,21)$ appear. $x_{lc}$ ranges between
$-1291$ and $-3767$, but beyond these $x_{lc}$ no cycle down to
$-2\cdot10^7$ appears.

\paragraph{Example $k=5$}
Five positive limit cycles were found for $k=5$,
one with $(a,b,D,N_{cyc})=(3,1,5,1)$,
two with $(a,b,D,N_{cyc})=(5,3,5,2)$, and
two with $(a,b)=(27,17)$.

\paragraph{Example $k=13$}
Only positive limit cycles are found for $k=13$:
one limit cycle with $(a,b,D,N_{cyc})=(4,1,13,1)$,
one with $(a,b)=(24,15)$ and
seven with $(a,b,D,N_{cyc})=(8,5,13,7)$.

\paragraph{Example $k=5137$}
The only limit cycles for $|x_{lc}|$ in the interval
$0.5\cdot10^7...2\cdot10^7$ were found for
$k=5137$. There is a large number of limit cycles with $(a,b)=(84,53)$
and $x_{lc}$ ranging from $-3056657$ down to $-17904629$.

\paragraph{Periodic lattice of $(a,b)$}
If $2^{a_1}\bmod k = 3^{b_1}\bmod k$ and
$2^{a_2}\bmod k = 3^{b_2}\bmod k$, then 
\be
2^{n_1a_1+n_2a_2}\bmod k = 3^{n_1b_1+n_2b_2}\bmod k
\ee
with $n_1,n_2\in\cZ$. Thus the appearing $(a,b)$ lie in
a region of a periodic lattice,
\be
(a,b)=n_1(a_1,b_1)+n_2(a_2,b_2).
\ee
$k=9215=5\cdot19\cdot97$ and $k=1843=19\cdot97$ have the same
$(a_1,b_1)=(27,9)$, $(a_2,b_2)=(17,11)$.
The corresponding $(a,b)$ and number of limit cycles $N_{a,b}$
are listed in table \ref{tb6}.\\
$k=485=5\cdot97$ and $k=97$ have $(a_1,b_1)=(9,3)$, $(a_2,b_2)=(17,11)$.\\
$k=95=5\cdot19$ and $k=19$ have $(a_1,b_1)=(11,5)$, $(a_2,b_2)=(17,11)$.
\begin{table}
\noindent
\begin{tabular}{rrrrr}
$(a,b)$ & $n_1$ & $n_2$ & $N^{(k)}_{a,b}$ & $N^{(k')}_{a,b}$
 \\ \hline
(17,11) & 0 & 1 & 92 & 22 \\
(34,22) & 0 & 2 & 14 & 7 \\
(51,33) & 0 & 3 & 6 \\
(68,44) & 0 & 4 & 2 \\
(27,9) & 1 & 0 & 11 & 2 \\
(44,20) & 1 & 1 & 18 & 4 \\
(61,31) & 1 & 2 & 7 & 2\\
(78,42) & 1 & 3 & 3 & 2 \\
\end{tabular}
\hspace{5mm}
\begin{tabular}{rrrrr}
$(a,b)$ & $n_1$ & $n_2$ & $N^{(k)}_{a,b}$ & $N^{(k')}_{a,b}$
 \\ \hline
(95,53) & 1 & 4 & 8 & 2 \\
(112,64) & 1 & 5 & 3 & 2 \\
(71,29) & 2 & 1 & 1 \\
(88,40) & 2 & 2 & 1 \\
(105,51) & 2 & 3 & 3 & 1 \\
(122,62) & 2 & 4 & 1 \\
(139,73) & 2 & 5 & 1 \\
(207,117) & 2 & 9 & & 1
\end{tabular}
\caption{Number $N_{a,b}$ of cycles for $k=9215$ and $k'=1843$
with $(a_1,b_1)=(27,9)$, $(a_2,b_2)=(17,11)$.} \label{tb6}
\end{table}

\section{Tao's recursion formula and Oscillator function}

In a recent paper \cite{Tao} Terrence Tao showed

Theorem 1.3 (Almost all Collatz orbits attain almost bounded
values) {\it Let $f:\cN+1\rightarrow \cR$ be any function
with $\lim_{N\rightarrow\infty} f(N)\rightarrow +\infty.$
Then one has $\Col_{min}(N)<f(N)$ for almost all $N\in\cN+1$.}
$\Col_{min}(N)$ denotes the minimal element of the Collatz orbit.

Here some aspects of his {\it
subsection 1.4 Fine-scale mixing of Syracuse random variables}
are considered for the more general case $3x+k$ with $k\in\cD$
in comparison to the mapping with $3x+1$.

\paragraph{Probability distribution} We define the probability
\be
P_{k,n}(y) = \lim_{N\rightarrow\infty} \frac 1N \sum_{z=1}^{N}
1(\Sir_k^n(2z+1) = y\bmod 3^n),
\ee
where $1$ is the indicator function, 1(true)=1, 1(false)=0.
We derive the recurrence relation: Let
\be
x_{n+1} = 3^{n+1} N + y, \quad x_n= 3^n N'+y'. \label{Recxk}
\ee
The iteration
\be
x_{n+1}=\Sir_k(x_n)=\frac{3x_n+k}{2^a}
\ee
yields
\be
N' = 2^a N, \quad y'=\frac{2^a y-k}3.
\ee
Hence $x_n$ is obtained from $x_{n+1}$ with probability $2^{-a}$.
Summation over all $a$, which yield integer $y'$, gives
\be
P_{k,n+1}(y) = \sum_{a\ge1:(2^ay-k)\bmod 3=0} 2^{-a}
P_{k,n}\left(\frac{2^a y-k}3\right). \label{RecP}
\ee
For $k=1$ this is the recursive formula for Syracuse random
variables Lemma 1.12 by Tao.
Tao's $\cP(\Syr...)$ is related to our $P$ by
\be
\cP(\Syr(\cZ/3^n\cZ)=x) = P_{1,n}(x).
\ee
and his recurrence relation is eq. (\ref{RecP}) for $k=1$.
Eq. (\ref{Recxk}) may be written
\be
k^{-1} x_{n+1} = \frac{3k^{-1}x_n+1}{2^a} = \Sir_1(k^{-1}x_n).
\label{Recxk2}
\ee
The set $k\in\cZ/3^n\cZ:3\nmid k$ forms a group under
multiplication. Thus the inverse $k^{-1}$ modulo $3^n$
for $k$ not divisible by 3 exists. For example one has
\be
1\cdot1 = 2\cdot5 = 4\cdot7 = 8\cdot8 = 1 \bmod 3^2.
\ee
Eq. (\ref{Recxk2}) yields
\be
P_{k,n}(x) = P_{1,n}(k^{-1}x), \quad\succ\quad
P_{k,n}(kx) = P_{1,n}(x). \label{Pk1n}
\ee
Thus the set of the probabilities $P_{k,n}$ and $P_{1,n}$ 
is the same but for different $x$.
The probabilities for $n=1,2,3$ are given in table \ref{tb3}
for $k=1$ and $k=5$. In subsequent lines $x$ is multiplied by 
2 modulo $3^n$.

\begin{table}

\begin{tabular}{rrr}
  x &  x & $P_{k,1}(x)$ \\
$k=1$ & $k=5$ & \\ \hline
  1 &  2 &  1/3 \\
  2 &  1 &  2/3 \\ \hline
 \rule{0mm}{5mm} & & \\
  x &  x & $P_{k,2}(x)$ \\ 
$k=1$ & $k=5$ & \\ \hline
  7 &  8 &  2/63 \\
  5 &  7 &  4/63 \\
  1 &  5 &  8/63 \\
  2 &  1 & 16/63 \\ \hline
  4 &  2 & 11/63 \\
  8 &  4 & 22/63 \\ \hline
\end{tabular}
\hspace{5mm}
\begin{tabular}{rrr}
  x &  x & $(2^{18}-1) P_{k,3}(x)$ \\ 
$k=1$ & $k=5$ & \\ \hline
  7 &  8 &  2408 \\
 14 & 16 &  4816 \\
  1 &  5 &  9632 \\
  2 & 10 & 19264 \\ \hline
  4 & 20 &  5240 \\
  8 & 13 & 10480 \\ \hline
 16 & 26 &  4316 \\
  5 & 25 &  8632 \\
 10 & 23 & 17264 \\
 20 & 19 & 34528 \\ \hline
 13 & 11 & 23285 \\
 26 & 22 & 46570 \\ \hline
 25 & 17 &  1598 \\
 23 &  7 &  3196 \\
 19 & 14 &  6392 \\
 11 &  1 & 12784 \\ \hline
 22 &  2 & 17246 \\
 17 &  4 & 34492 \\ \hline
\end{tabular}

\caption{The probabilities $P_{k,n}(x)$ for $n=1,2,3$.
($2^{18}-1=63\cdot4161$)} \label{tb3}

\end{table}

More generally we consider the probability of $\Sir^l$ for $x\bmod 3^n$,
\be
Q_{k,l,n}(x) = \lim_{N\rightarrow\infty} \frac1N
\sum_{z=1}^N 1(\Sir_k^l(2z+1) = x\bmod3^n).
\ee

\paragraph{Case $l>n$} In this case the probability is given by
\be
Q_{k,l,n}(x) = \sum_{z=0}^{3^{l-n}-1} P_{k,l}(x+3^n z). \label{Qkln}
\ee
Then one obtains
\bea
Q_{k,l+1,n}(x) &=& \sum_{z=0}^{3^{l+1-n}-1} 
\sum_{a\ge1:(2^ax-k)\bmod3=0}
2^{-a} P_{k,l}\left(\frac{2^a(x+3^n z)-k}3\right) \\
&=& \sum_{a\ge1:(2^ax-k)\bmod 3=0} 2^{-a} \sum_{z=0}^{3^{l+1-n}-1}
P_{k,l}\left(\frac{2^ax-k}3 +3^{n-1}2^a z\right). \nonumber
\eea
The set of elements $2^az$ from $z=0$ to $3^{l+1-n}-1$ agrees
modulo $3^{l+1-n}$ with the set of elements $z$.
Hence one obtains
\bea
\sum_{z=0}^{3^{l+1-n}-1} P_{k,l}\left(\frac{2^ax-k}3 +3^{n-1}2^a z\right)
&=& \sum_{z=0}^{3^{l+1-n}-1} P_{k,l}\left(\frac{2^ax-k}3 +3^{n-1} z\right) \nn
&=& Q_{k,l,n-1}\left(\frac{2^ax-k}3\right)
\eea
and thus the recurrence relation
\be
Q_{k,l+1,n}(x) = \sum_{a\ge1:(2^ax-k)\bmod3=0} 2^{-a} Q_{k,l,n-1}
\left(\frac{2^ax-k}3\right). \label{Qkln1}
\ee
Thus $Q$ obeys the same recurrence relation as $P$.
The definition (\ref{Qkln}) yields
\be
Q_{k,l-n,0}(x) = \sum_{z=0}^{3^{l-n}-1} P_{k,l-n}(x+z) = 1 = P_{k,0}(0),
\ee
since the summation is over all $z$ modulo $3^{n-l}$. From
this initial condition and the recurrence relation (\ref{Qkln1})
one obtains
\be
Q_{k,l,n}(x) = P_{k,n}(x).
\ee
Thus the application of the mapping $\Sir$ beyond $\Sir^n$
does not change the probability distribution modulo $3^n$.

\paragraph{Case $l<n$}. In this case the probability $Q_{k,l,n}(x)$
is invariant for $x$ modulo $3^l$. Thus
\be
Q_{k,l,n}(x) = 3^{l-n} P_{k,l}(x).
\ee
In total
\be
Q_{k,l,n}(x) = \left\{ \begin{array}{rr} P_{k,n}(x) & l\ge n \\
3^{l-n} P_{k,l}(x) & l\le n \end{array} \right. .
\ee

\paragraph{Oscillation} Tao\cite{Tao} introduces the oscillation in 
proposition 1.14 (Fine scale mixing of n-Syracuse offsets)
\bea
\Osc_{m,n}(c_Y)_{Y\in\cZ/3^n\cZ} = \sum_{Y\in\cZ/3^n\cZ}
\left| c_Y-3^{m-n}\sum_{Y'\in\cZ/3^n\cZ:Y'=Y\bmod 3^m} c_{Y'}\right|.
\eea
with $c_Y=P_{k,n}(Y)$ with $k=1$. The oscillation is for
general $k$
\be
\Osc_{m,n} = \sum_{Y\in\cZ/3^n\cZ}
| P_{k,n}(Y) - 3^{m-n} P_{k,m}(Y) |.
\ee
With (\ref{Pk1n}) we have $P_{k,n}(Y)=P_{1,n}(k^{-1}Y)$,
similarly for $P_{k,m}$. 
$Y$ runs through all $3^n$ elements of $\cZ/3^n\cZ$.
Since the mapping between $Y$ and $k^{-1}Y$ is bijective,
also $k^{-1}Y$ runs through all $3^n$ elements of $\cZ/3^n\cZ$.
Thus the oscillation $\Osc_{m,n}$ does not depend on $k\in\cD$.

From table \ref{tb3} one obtains
\be
\Osc_{1,2} = 10/21, \quad \Osc_{1,3} = 2, \quad \Osc_{2,3} = 2/3.
\ee
The properties of the oscillation $\Osc_{m,n}$ are essential
in Tao's proof of his theorem 1.3: {\it Almost all
Collatz orbits attain almost bounded values}.


\begin{thebibliography}{9}

\bibitem{BelMig06a} Edward G. Belaga, Maurice Mignotte, {\it Walking cautiously
into the Collatz Wilderness: Algorithmically, Number-Theoretically,
randomly}, (2006) hal.archives-ouvertes.fr/hal-00129726

\bibitem{BelMig06b} Edward G. Belaga, Maurice Mignotte, 
{\it The Collatz Problem and its Generalizations: Experimental
Data. Table 1. Primitive Cycles of $(3n+d)$-mappings},
(2006) hal.archives-ouvertes.fr/hal-00129727

\bibitem{BelMig12} Edward G. Belaga, Maurice Mignotte,
{\it Embedding the $3x+1$ conjecture in a $3x+d$ context},
Exp. Math. 7 (2012) 145-151

\bibitem{Bre02} Barry Brent, {\it $3x+1$ dynamics on rationals
with fixed denominator}, arXiv:math/0204170v1

\bibitem{Lag90} Jeffrey C. Lagarias, {\it The set of rational
cycles for the $3x+1$ problem}m Acta Arithmetica 56 (1990)
33-53

\bibitem{Lag03} Jeffrey C. Lagarias, {\it The $3x+1$ Problem:
An annotated bibliography (1963-1999) (sorted by author)}
arXiv:math/0309224

\bibitem{Lag06} Jeffrey C. Lagarias, {\it The $3x+1$ Problem:
An annotated bibliography, II (2000-2009)}
arXiv:math/0608208

\bibitem{Tao} Terrence Tao, {\it Almost all orbits of the Collatz map
attain almost bounded values} (2019) arXiv:1909.03562v2

\bibitem{Weglist} Franz Wegner, {\it Supplement: List of limit
cycles for the Collatz problem $3x+k$}

\end{thebibliography}
\end{document}


\title{Supplement: List of limit cycles for the Collatz problem $3x+k$}

\author{Franz Wegner \\
Institut f\"ur Theoretische Physik \\
Universit\"at Heidelberg \\
D-69120 Heidelberg Germany}

\maketitle

\begin{abstract}
In this supplement to "The Collatz problem generalized to $3x+k$"
a list of limit cycles is given for $k$ not divisible by 2 and 3.
Both positive and negative limit cycles are presented up to $k=9997$
starting with $x_0=-2\cdot10^7...+2\cdot10^7$.
\end{abstract}

\section*{Introduction}

The Collatz sequence starts with $x_0$ and defines a sequence
$x_1,x_2,...$ by the rule
$$
x_{i+1} = {\rm Col}_k(x_i), \quad
{\rm Col}_k(x) := \left\{ \begin{array}{cc} 3x+k, & x \mbox{ odd}, \\
x/2 & x \mbox{ even} \end{array} \right.
$$
For reasons explained in the main paper we consider positive $k$ not
divisible by 2 and 3 and start from $x_0$ not divisible by 2 and 3
and prime to $k$. The calculation is restricted to positive $k$
less than 10000 and to $x_0$ in the interval
$-2\cdot10^7...+2\cdot10^7$.
In all these cases the sequence runs into a limit case. This means
that for each $x_0$ there exist constants $n_0$ and $c$ so that
$x_{n+c}=x_n$ for $n\ge n_0$. (The smallest positive $c$ is chosen).
The sequence from $x_n$ to $x_{n+c}$
contains $a$ steps $x/2$ and $b$ steps $3x+k$, $c=a+b$.

For each $k$ the numbers $x_0$, $x_{lc}$, $a$, and $b$ are listed
for each limit cycle.
There are limit cycles with positive $x$ and those with
negative $x$. I call them positive and negative limit cycles.
For the positive limit cycles I list the smallest $x$ in
the cycle and denote it by $x_{lc}$. For negative limit cycles
$x_{lc}$ is that with the smallest $|x|$ in the cycle.
Also $x_0$ is given. For positive
limit cycles it is the smallest positive $x_0$ from which the
cycle is reached, for negative cycles it is the smallest
$|x_0|$ from which it is reached.

\section{Limit cycles: k from 1 to 997}

\noindent
\begin{tabular}{rrrrr}
 $k$ & $x_0$ & $x_{lc}$ & $a$ & $b$ \\
\hline
  1 &     1 &     1 &    2 &    1 \\
  1 &    -1 &    -1 &    1 &    1 \\
  1 &    -5 &    -5 &    3 &    2 \\
  1 &   -17 &   -17 &   11 &    7 \\ \hline
  5 &     1 &     1 &    3 &    1 \\
  5 &     7 &    19 &    5 &    3 \\
  5 &    23 &    23 &    5 &    3 \\
  5 &   187 &   187 &   27 &   17 \\
  5 &   259 &   347 &   27 &   17 \\ \hline
  7 &     1 &     5 &    4 &    2 \\ \hline
 11 &     1 &     1 &    6 &    2 \\
 11 &     5 &    13 &   14 &    8 \\
 11 &   -19 &   -19 &    4 &    3 \\ \hline
 13 &     1 &     1 &    4 &    1 \\
 13 &    19 &   131 &   24 &   15 \\
 13 &   155 &   211 &    8 &    5 \\
 13 &   163 &   251 &    8 &    5 \\
 13 &   187 &   287 &    8 &    5 \\
 13 &   191 &   259 &    8 &    5 \\
 13 &   227 &   227 &    8 &    5 \\
 13 &   283 &   283 &    8 &    5 \\
 13 &   319 &   319 &    8 &    5 \\ \hline
 17 &     1 &     1 &    7 &    2 \\
 17 &    11 &    23 &   31 &   18 \\
 17 &   -41 &   -65 &    6 &    4 \\
 17 &   -73 &   -73 &    6 &    4 \\ \hline
 19 &     1 &     5 &   11 &    5 \\
 19 &   -67 &  -115 &   17 &   11 \\ \hline
 23 &     1 &    41 &   43 &   26 \\
 23 &     5 &     5 &    5 &    2 \\
 23 &     7 &     7 &    5 &    2 \\
 23 & -2263 & -2263 &   19 &   12 \\
 23 & -2359 & -2359 &   19 &   12 \\
 23 & -2495 & -2963 &   19 &   12 \\
 23 & -2503 & -3743 &   19 &   12 \\
 23 & -2567 & -3415 &   19 &   12 \\
 23 & -2743 & -2743 &   19 &   12 \\
 23 & -3091 & -3091 &   19 &   12 \\
 23 & -4819 & -4819 &   19 &   12 \\ \hline
\end{tabular}
\hspace{5mm}
\noindent\begin{tabular}{rrrrr}
 $k$ & $x_0$ & $x_{lc}$ & $a$ & $b$ \\
\hline
 25 &     1 &     7 &   16 &    8 \\
 25 &    17 &    17 &    8 &    4 \\
 25 &   -89 &  -113 &   34 &   22 \\ \hline
 29 &     1 &     1 &    5 &    1 \\
 29 &     5 &    11 &   17 &    9 \\
 29 &  2531 &  3811 &   65 &   41 \\
 29 &  7055 &  7055 &   65 &   41 \\
 29 &   -85 &  -109 &   12 &    8 \\ \hline
 31 &     1 &    13 &   23 &   12 \\
 31 &   -79 &   -95 &    9 &    6 \\ \hline
 35 &     1 &    13 &    8 &    4 \\
 35 &    11 &    17 &    8 &    4 \\ \hline
 37 &     1 &    19 &    6 &    3 \\
 37 &     7 &    29 &    6 &    3 \\
 37 &    23 &    23 &    6 &    3 \\ \hline
 41 &     1 &     1 &   20 &    8 \\ \hline
 43 &     1 &     1 &   11 &    3 \\ \hline
 47 &     1 &    25 &   28 &   16 \\
 47 &     5 &     5 &   18 &    7 \\
 47 &    37 &    65 &    7 &    4 \\
 47 &    41 &    85 &    7 &    4 \\
 47 &    53 &    89 &    7 &    4 \\
 47 &    73 &    73 &    7 &    4 \\
 47 &   101 &   101 &    7 &    4 \\ \hline
 49 &     1 &    25 &   38 &   22 \\
 49 &   -65 &   -65 &    5 &    4 \\ \hline
 53 &     1 &   103 &   29 &   17 \\ \hline
 55 &     1 &     1 &   12 &    4 \\
 55 &     7 &     7 &    6 &    2 \\
 55 &    41 &    41 &   28 &   16 \\ \hline
 59 &     1 &     1 &   28 &   11 \\
 59 &     5 &   133 &   10 &    6 \\
 59 &    61 &   181 &   10 &    6 \\
 59 &   125 &   217 &   10 &    6 \\
 59 &   149 &   149 &   10 &    6 \\
 59 &   169 &   185 &   10 &    6 \\
 59 &   221 &   221 &   10 &    6 \\ \hline
 61 &     1 &     1 &    6 &    1 \\
 61 &   175 &   235 &   66 &   41 \\ \hline
\end{tabular}
\newpage
\noindent\begin{tabular}{rrrrr}
 $k$ & $x_0$ & $x_{lc}$ & $a$ & $b$ \\
\hline
 65 &     1 &    19 &   24 &   12 \\ \hline
 67 &     1 &    17 &   30 &   16 \\ \hline
 71 &     1 &    29 &   10 &    5 \\
 71 &     7 &    31 &   10 &    5 \\
 71 &   893 &  4409 &   27 &   17 \\
 71 &  1073 &  2809 &   27 &   17 \\
 71 &  2585 &  2585 &   27 &   17 \\
 71 &  2633 &  3985 &   27 &   17 \\
 71 &  4121 &  4121 &   27 &   17 \\ \hline
 73 &     1 &    19 &   60 &   32 \\
 73 &     5 &     5 &   12 &    4 \\
 73 &     7 &    47 &   15 &    8 \\ \hline
 77 &     1 &     1 &   38 &   16 \\
 77 &  -397 & -1153 &   22 &   14 \\
 77 &  -685 &  -989 &   22 &   14 \\
 77 & -1165 & -1165 &   22 &   14 \\ \hline
 79 &     1 &     1 &   44 &   20 \\
 79 &     7 &     7 &   44 &   20 \\
 79 &   193 &   265 &   23 &   14 \\
 79 &   233 &   233 &   23 &   14 \\ \hline
 83 &     1 &   109 &   12 &    7 \\
 83 &    17 &   157 &   12 &    7 \\
 83 &    59 &    65 &   24 &   14 \\ \hline
 85 &     1 &     7 &  100 &   56 \\
 85 &  -293 &  -293 &   12 &    8 \\ \hline
 89 &     1 &    17 &   17 &    8 \\ \hline
 91 &     1 &     1 &   48 &   24 \\
 91 &    25 &    25 &   12 &    6 \\
 91 &    59 &    59 &   12 &    6 \\ \hline
 95 &     1 &     1 &   72 &   36 \\
 95 &     7 &    23 &   11 &    5 \\
 95 &    17 &    17 &   11 &    5 \\
 95 &  -511 &  -511 &   17 &   11 \\ \hline
 97 &     1 &     1 &   18 &    6 \\
 97 &     5 &    13 &    9 &    3 \\
 97 &  -449 &  -625 &   17 &   11 \\
 97 &  -505 &  -505 &   17 &   11 \\ \hline
101 &     1 &    11 &   14 &    6 \\
101 &     5 &    29 &    7 &    3 \\
101 &     7 &     7 &   14 &    6 \\
101 &    17 &    19 &    7 &    3 \\
101 &    23 &    23 &    7 &    3 \\
101 &    31 &    31 &    7 &    3 \\
\end{tabular}
\hspace{5mm}
\noindent\begin{tabular}{rrrrr}
 $k$ & $x_0$ & $x_{lc}$ & $a$ & $b$ \\
\hline
101 &    37 &    37 &    7 &    3 \\ \hline
103 &     1 &    23 &   45 &   22 \\
103 &     5 &     5 &   21 &    8 \\ \hline
107 &     1 &     1 &  106 &   53 \\ \hline
109 &     1 &    19 &   52 &   30 \\ \hline
113 &     1 &     1 &   29 &   12 \\
113 &  -769 &  -769 &   69 &   44 \\ \hline
115 &     1 &    13 &   27 &   13 \\
115 &     7 &    17 &   34 &   18 \\
115 &  -163 &  -211 &    7 &    5 \\
115 &  -227 &  -227 &    7 &    5 \\
115 &  -251 &  -251 &    7 &    5 \\ \hline
119 &     1 &     1 &   75 &   42 \\
119 &     5 &     5 &    7 &    2 \\
119 &    11 &    11 &    7 &    2 \\
119 &    13 &    23 &   39 &   18 \\
119 &   125 &   125 &   25 &   14 \\
119 &  -247 &  -311 &   18 &   12 \\ \hline
121 &     1 &     5 &   44 &   18 \\
121 &    19 &    19 &   88 &   46 \\
121 & -1085 & -1085 &   44 &   28 \\
121 & -1105 & -2089 &   22 &   14 \\ \hline
125 &     1 &     1 &    7 &    1 \\
125 &     7 &    47 &   33 &   19 \\
125 &   143 &   143 &   33 &   19 \\
125 &   431 &   899 &  118 &   74 \\ \hline
127 &     1 &     1 &   37 &   18 \\
127 &     5 &    41 &   37 &   18 \\ \hline
131 &     1 &    13 &   26 &   13 \\
131 &     7 &    17 &   26 &   13 \\
131 &    23 &    23 &   52 &   26 \\
131 & -7471 & -7471 &   38 &   24 \\
131 & -8003 & -15875 &   38 &   24 \\
131 & -10883 & -10883 &   38 &   24 \\
131 & -12107 & -16099 &   38 &   24 \\
131 & -13963 & -13963 &   38 &   24 \\ \hline
133 &     1 &    11 &   36 &   18 \\
133 &    31 &    59 &   42 &   24 \\ \hline
137 &     1 &     1 &   14 &    4 \\
137 &     5 &    41 &   72 &   40 \\
137 &    95 &   503 &   29 &   18 \\
137 &   415 &   967 &   29 &   18 \\
137 &   743 &   743 &   29 &   18 \\ \hline
\end{tabular}
\newpage
\noindent\begin{tabular}{rrrrr}
 $k$ & $x_0$ & $x_{lc}$ & $a$ & $b$ \\
\hline
139 &     1 &    11 &  136 &   74 \\
139 &  -331 & -1291 &   22 &   14 \\
139 &  -587 & -1931 &   22 &   14 \\
139 &  -779 & -2059 &   11 &    7 \\
139 &  -827 & -2203 &   22 &   14 \\
139 & -1043 & -1495 &   33 &   21 \\
139 & -1195 & -2251 &   11 &    7 \\
139 & -1339 & -2539 &   11 &    7 \\
139 & -1355 & -1355 &   22 &   14 \\
139 & -1367 & -1367 &   33 &   21 \\
139 & -1555 & -2971 &   11 &    7 \\
139 & -1835 & -2683 &   11 &    7 \\
139 & -1915 & -2507 &   11 &    7 \\
139 & -1979 & -2899 &   11 &    7 \\
139 & -2123 & -2123 &   11 &    7 \\
139 & -2131 & -2795 &   11 &    7 \\
139 & -2171 & -3187 &   11 &    7 \\
139 & -2195 & -3223 &   11 &    7 \\
139 & -2219 & -2219 &   11 &    7 \\
139 & -2347 & -2347 &   11 &    7 \\
139 & -2387 & -3511 &   11 &    7 \\
139 & -2419 & -2903 &   11 &    7 \\
139 & -2455 & -3227 &   11 &    7 \\
139 & -2491 & -2491 &   11 &    7 \\
139 & -2519 & -2747 &   11 &    7 \\
139 & -2579 & -2579 &   11 &    7 \\
139 & -2603 & -2603 &   11 &    7 \\
139 & -2707 & -2707 &   11 &    7 \\
139 & -2711 & -2963 &   11 &    7 \\
139 & -2743 & -3287 &   11 &    7 \\
139 & -2939 & -2939 &   11 &    7 \\
139 & -3031 & -3031 &   11 &    7 \\
139 & -3155 & -3155 &   11 &    7 \\
139 & -3443 & -3443 &   11 &    7 \\
139 & -3479 & -3479 &   11 &    7 \\
139 & -3767 & -3767 &   11 &    7 \\ \hline
143 &     1 &     7 &  140 &   80 \\
143 &     5 &    29 &   16 &    7 \\
143 &    17 &    17 &   16 &    7 \\ \hline
145 &     1 &     1 &   10 &    2 \\
145 &     7 &    23 &   22 &   10 \\
145 &    17 &    47 &   34 &   18 \\
145 &  -209 &  -241 &   48 &   32 \\
\end{tabular}
\hspace{5mm}
\noindent\begin{tabular}{rrrrr}
 $k$ & $x_0$ & $x_{lc}$ & $a$ & $b$ \\
\hline
145 &  -617 &  -617 &   12 &    8 \\ \hline
149 &     1 &    19 &   59 &   33 \\
149 &   395 &   667 &   42 &   26 \\ \hline
151 &     1 &    41 &   72 &   40 \\
151 &     5 &   113 &   18 &   10 \\ \hline
155 &     1 &     1 &   42 &   18 \\ \hline
157 &     1 &    13 &   28 &   12 \\ \hline
161 &     1 &    25 &   44 &   22 \\
161 &     5 &    11 &   15 &    6 \\
161 &    19 &    19 &   15 &    6 \\
161 &    65 &    79 &   29 &   16 \\ \hline
163 &     1 &    11 &   89 &   49 \\
163 &     7 &    17 &   72 &   36 \\ \hline
167 &     1 &   127 &   17 &    9 \\
167 &     5 &    13 &   47 &   20 \\
167 &    71 &    95 &   17 &    9 \\ \hline
169 &     1 &    11 &   24 &    9 \\
169 &     7 &    17 &   44 &   23 \\ \hline
173 &     1 &     7 &   34 &   14 \\
173 &     5 &    37 &   17 &    7 \\ \hline
175 &     1 &    89 &    8 &    4 \\
175 &    11 &   103 &    8 &    4 \\
175 &    17 &    29 &   24 &   12 \\
175 &    37 &   143 &    8 &    4 \\
175 &    53 &    53 &   16 &    8 \\
175 &    73 &    73 &    8 &    4 \\
175 &   101 &   101 &    8 &    4 \\ \hline
179 &     1 &     5 &  110 &   62 \\
179 &  -211 &  -211 &    6 &    5 \\ \hline
181 &     1 &    23 &   20 &   10 \\
181 &     5 &    11 &   80 &   40 \\
181 &    13 &    55 &   20 &   10 \\
181 & -3305 & -3605 &   68 &   43 \\
181 & -3845 & -5677 &   68 &   43 \\ \hline
185 &     1 &     1 &   72 &   36 \\
185 & -1273 & -1273 &   44 &   28 \\ \hline
187 &     1 &    43 &   26 &   12 \\
187 &     5 &     5 &   26 &   12 \\ \hline
191 &     1 &   961 &   40 &   25 \\
191 &     7 &    47 &   46 &   24 \\
191 &  1153 &  2405 &   40 &   25 \\ \hline
193 &     1 &     1 &   66 &   29 \\
193 &     5 &     5 &   18 &    5 \\
\end{tabular}
\newpage
\noindent\begin{tabular}{rrrrr}
 $k$ & $x_0$ & $x_{lc}$ & $a$ & $b$ \\
\hline
193 &    11 &    91 &   84 &   50 \\
193 &    31 &    31 &   24 &   12 \\
193 & -2393 & -5617 &   30 &   19 \\
193 & -3329 & -4897 &   30 &   19 \\
193 & -3841 & -3841 &   30 &   19 \\
193 & -4561 & -6017 &   30 &   19 \\
193 & -4817 & -4817 &   30 &   19 \\
193 & -6053 & -6053 &   30 &   19 \\ \hline
197 &     1 &     5 &  141 &   69 \\ \hline
199 &     1 &    13 &   15 &    6 \\
199 &     5 &    47 &   15 &    6 \\ \hline
203 &     1 &     1 &   10 &    2 \\
203 &     5 &    43 &  102 &   54 \\ \hline
205 &     1 &    43 &   80 &   40 \\ \hline
209 &     1 &     1 &   68 &   26 \\
209 &    13 &   235 &   14 &    8 \\
209 &    47 &   175 &   14 &    8 \\
209 &   203 &   359 &   14 &    8 \\ \hline
211 &     1 &    37 &   88 &   46 \\ \hline
215 &     1 &     1 &   38 &   18 \\
215 &    11 &    31 &   19 &    9 \\
215 &    13 &    47 &   19 &    9 \\ \hline
217 &     1 &    17 &   51 &   24 \\
217 &     5 &    59 &   78 &   42 \\
217 &  -281 &  -313 &   54 &   36 \\
217 &  -569 &  -745 &    9 &    6 \\
217 &  -617 &  -817 &    9 &    6 \\
217 &  -625 &  -761 &    9 &    6 \\
217 &  -689 &  -905 &    9 &    6 \\
217 &  -697 &  -697 &    9 &    6 \\
217 &  -809 &  -809 &    9 &    6 \\
217 &  -881 &  -881 &    9 &    6 \\
217 &  -977 &  -977 &    9 &    6 \\ \hline
221 &     1 &     5 &   20 &    8 \\ \hline
223 &     1 &    71 &   47 &   24 \\
223 &    19 &    35 &   47 &   24 \\ \hline
227 &     1 &     5 &   66 &   26 \\ \hline
229 &     1 &    29 &    8 &    3 \\
229 &     5 &     7 &   16 &    6 \\
229 &    19 &    19 &    8 &    3 \\
229 &    23 &    23 &    8 &    3 \\
229 &    31 &    31 &    8 &    3 \\
229 &    37 &    37 &    8 &    3 \\
\end{tabular}
\hspace{5mm}
\noindent\begin{tabular}{rrrrr}
 $k$ & $x_0$ & $x_{lc}$ & $a$ & $b$ \\
\hline
229 &    47 &    47 &    8 &    3 \\
229 &    53 &    53 &    8 &    3 \\ \hline
233 &     1 &     5 &  107 &   48 \\
233 &    23 &   919 &   13 &    8 \\
233 &   191 &  1015 &   26 &   16 \\
233 &   209 &  1111 &   13 &    8 \\
233 &   407 &   727 &   26 &   16 \\
233 &   647 &   647 &   26 &   16 \\
233 &   695 &  1159 &   13 &    8 \\
233 &   787 &  1663 &   13 &    8 \\
233 &   823 &  1175 &   13 &    8 \\
233 &   959 &  1555 &   13 &    8 \\
233 &   967 &  1367 &   13 &    8 \\
233 &   983 &   983 &   13 &    8 \\
233 &  1079 &  1079 &   13 &    8 \\
233 &  1183 &  1655 &   13 &    8 \\
233 &  1223 &  1223 &   13 &    8 \\
233 &  1235 &  1535 &   13 &    8 \\
233 &  1439 &  1439 &   13 &    8 \\
233 &  1651 &  1651 &   13 &    8 \\
233 &  1907 &  1907 &   13 &    8 \\ \hline
235 &     1 &    17 &   25 &   11 \\
235 &    11 &    11 &   25 &   11 \\
235 &    13 &    19 &   25 &   11 \\ \hline
239 &     1 &     1 &   58 &   26 \\
239 &     7 &    13 &   29 &   13 \\
239 &    19 &    19 &   29 &   13 \\
239 &    55 &   113 &   31 &   18 \\ \hline
241 &     1 &    59 &   54 &   30 \\
241 &     7 &     7 &   80 &   40 \\ \hline
245 &     1 &    43 &  114 &   66 \\
245 &   167 &   167 &   38 &   22 \\ \hline
247 &     1 &     1 &   72 &   36 \\
247 &     5 &     5 &    8 &    2 \\
247 &     7 &     7 &    8 &    2 \\
247 &    11 &    11 &    8 &    2 \\
247 &    17 &   145 &   28 &   16 \\
247 &   191 &   337 &   28 &   16 \\ \hline
251 &     1 &    91 &   74 &   40 \\
251 &   917 &  1445 &   80 &   50 \\ \hline
253 &     1 &     1 &    8 &    1 \\
253 &     5 &    17 &   42 &   19 \\
253 &    13 &    13 &  162 &   89 \\ \hline
\end{tabular}
\newpage
\noindent\begin{tabular}{rrrrr}
 $k$ & $x_0$ & $x_{lc}$ & $a$ & $b$ \\
\hline
257 &     1 &     1 &   43 &   16 \\
257 &   367 &   367 &   27 &   16 \\
257 &   635 &   635 &   27 &   16 \\ \hline
259 &     1 &    67 &   36 &   18 \\
259 &     5 &   107 &   12 &    6 \\
259 &    13 &   121 &   12 &    6 \\
259 &    17 &    59 &   24 &   12 \\
259 &    29 &    29 &   24 &   12 \\
259 &    85 &    85 &   12 &    6 \\
259 &   113 &   113 &   12 &    6 \\ \hline
263 &     1 &    17 &   40 &   21 \\
263 &     5 &     7 &  102 &   47 \\
263 &    11 &    37 &   40 &   21 \\ \hline
265 &     1 &    67 &   44 &   24 \\
265 &     7 &    47 &   30 &   14 \\
265 &    11 &    49 &   58 &   34 \\ \hline
269 &     1 &   331 &    9 &    5 \\
269 &     7 &    55 &   36 &   20 \\
269 &    13 &   179 &   18 &   10 \\
269 &    37 &   139 &   27 &   15 \\
269 &    59 &   223 &   18 &   10 \\
269 &    83 &   259 &    9 &    5 \\
269 &    91 &   211 &    9 &    5 \\
269 &   155 &   367 &    9 &    5 \\
269 &   175 &   323 &    9 &    5 \\
269 &   191 &   383 &    9 &    5 \\
269 &   227 &   227 &    9 &    5 \\
269 &   229 &   395 &    9 &    5 \\
269 &   245 &   251 &    9 &    5 \\
269 &   283 &   283 &    9 &    5 \\
269 &   287 &   287 &    9 &    5 \\
269 &   319 &   319 &    9 &    5 \\
269 &   347 &   347 &    9 &    5 \\
269 &   431 &   431 &    9 &    5 \\ \hline
271 &     1 &    17 &   45 &   20 \\ \hline
275 &     1 &     1 &  110 &   50 \\ \hline
277 &     1 &    13 &  120 &   63 \\ \hline
281 &     1 &     1 &   22 &    8 \\
281 &     5 &    37 &   11 &    4 \\
281 &    17 &    17 &   11 &    4 \\
281 &    19 &    19 &   11 &    4 \\
281 &    23 &    41 &   11 &    4 \\
281 &    25 &    25 &   11 &    4 \\
\end{tabular}
\hspace{5mm}
\noindent\begin{tabular}{rrrrr}
 $k$ & $x_0$ & $x_{lc}$ & $a$ & $b$ \\
\hline
281 &   145 &   287 &  106 &   64 \\ \hline
283 &     1 &    53 &   89 &   51 \\
283 &     5 &     5 &   16 &    6 \\
283 &   821 &   821 &   73 &   45 \\ \hline
287 &     1 &   169 &   25 &   14 \\
287 &     5 &     5 &   30 &   12 \\
287 &    25 &   121 &   25 &   14 \\
287 &    29 &    31 &   30 &   12 \\
287 &   127 &   229 &   25 &   14 \\
287 &   193 &   193 &   25 &   14 \\ \hline
289 &     1 &     1 &   33 &   14 \\
289 &     5 &    13 &   33 &   14 \\ \hline
293 &     1 &    11 &   41 &   17 \\ \hline
295 &     1 &   601 &   20 &   12 \\
295 &     7 &   737 &   20 &   12 \\
295 &    23 &   281 &   30 &   18 \\
295 &   233 &   761 &   10 &    6 \\
295 &   553 &   977 &   10 &    6 \\
295 &   593 &   889 &   10 &    6 \\
295 &   697 &   697 &   10 &    6 \\
295 &   701 &  1033 &   10 &    6 \\
295 &   721 &  1229 &   10 &    6 \\
295 &   733 &  1009 &   10 &    6 \\
295 &   809 &   809 &   10 &    6 \\
295 &   817 &   817 &   10 &    6 \\
295 &   829 &  1117 &   10 &    6 \\
295 &   881 &   881 &   10 &    6 \\
295 &   937 &   937 &   10 &    6 \\
295 &   989 &   989 &   10 &    6 \\
295 &  1213 &  1213 &   10 &    6 \\ \hline
299 &     1 &   197 &   36 &   21 \\
299 &     5 &     5 &   40 &   16 \\
299 &    11 &    19 &   20 &    8 \\ \hline
301 &     1 &    19 &   30 &   12 \\
301 &    61 &    67 &  108 &   60 \\ \hline
305 &     1 &    11 &   36 &   16 \\
305 &     7 &     7 &   36 &   16 \\
305 &    13 &    13 &   36 &   16 \\ \hline
307 &     1 &    17 &   21 &    9 \\
307 &     5 &     5 &   21 &    9 \\
307 &     7 &    41 &   21 &    9 \\
307 &    13 &   137 &   33 &   19 \\
307 &    35 &   205 &   33 &   19 \\
\end{tabular}
\newpage
\noindent\begin{tabular}{rrrrr}
 $k$ & $x_0$ & $x_{lc}$ & $a$ & $b$ \\
\hline
307 &   133 &   353 &   78 &   48 \\
307 &   293 &   293 &   33 &   19 \\ \hline
311 &     1 &    13 &   62 &   31 \\
311 &    11 &    29 &   62 &   31 \\
311 & -2359 & -44831 &   19 &   12 \\
311 & -4151 & -23351 &   19 &   12 \\
311 & -5683 & -26807 &   19 &   12 \\
311 & -6391 & -18551 &   38 &   24 \\
311 & -6967 & -24115 &   38 &   24 \\
311 & -7895 & -21791 &   38 &   24 \\
311 & -8615 & -15959 &   57 &   36 \\
311 & -8671 & -21239 &   38 &   24 \\
311 & -9011 & -29207 &   38 &   24 \\
311 & -9463 & -19799 &   57 &   36 \\
311 & -10663 & -31367 &   38 &   24 \\
311 & -11063 & -16439 &   38 &   24 \\
311 & -11191 & -32951 &   19 &   12 \\
311 & -13367 & -26423 &   38 &   24 \\
311 & -14303 & -28295 &   38 &   24 \\
311 & -14815 & -36791 &   19 &   12 \\
311 & -15319 & -22823 &   38 &   24 \\
311 & -15415 & -15415 &   38 &   24 \\
311 & -16823 & -33335 &   19 &   12 \\
311 & -17143 & -25559 &   38 &   24 \\
311 & -17831 & -35351 &   19 &   12 \\
311 & -17887 & -44567 &   19 &   12 \\
311 & -18487 & -27575 &   19 &   12 \\
311 & -18803 & -26327 &   19 &   12 \\
311 & -19255 & -19255 &   38 &   24 \\
311 & -19415 & -38519 &   19 &   12 \\
311 & -19507 & -32471 &   19 &   12 \\
311 & -19687 & -36823 &   19 &   12 \\
311 & -19703 & -26167 &   19 &   12 \\
311 & -20647 & -30815 &   38 &   24 \\
311 & -20983 & -20983 &   38 &   24 \\
311 & -21079 & -31463 &   38 &   24 \\
311 & -21719 & -28855 &   19 &   12 \\
311 & -22007 & -43703 &   19 &   12 \\
311 & -22951 & -34271 &   19 &   12 \\
311 & -24151 & -36071 &   19 &   12 \\
311 & -24743 & -49175 &   19 &   12 \\
311 & -24823 & -37079 &   19 &   12 \\
311 & -25139 & -41975 &   19 &   12 \\
\end{tabular}
\hspace{5mm}
\noindent\begin{tabular}{rrrrr}
 $k$ & $x_0$ & $x_{lc}$ & $a$ & $b$ \\
\hline
311 & -25399 & -25399 &   19 &   12 \\
311 & -25895 & -51479 &   19 &   12 \\
311 & -26551 & -26551 &   19 &   12 \\
311 & -26911 & -26911 &   38 &   24 \\
311 & -26975 & -53639 &   19 &   12 \\
311 & -27127 & -27127 &   19 &   12 \\
311 & -27443 & -45863 &   19 &   12 \\
311 & -27895 & -27895 &   19 &   12 \\
311 & -28127 & -37399 &   19 &   12 \\
311 & -28435 & -35575 &   19 &   12 \\
311 & -28919 & -34355 &   19 &   12 \\
311 & -29351 & -39031 &   19 &   12 \\
311 & -29623 & -29623 &   19 &   12 \\
311 & -29747 & -39559 &   19 &   12 \\
311 & -29911 & -29911 &   19 &   12 \\
311 & -29983 & -44819 &   19 &   12 \\
311 & -30259 & -33847 &   19 &   12 \\
311 & -30295 & -45287 &   19 &   12 \\
311 & -30431 & -40471 &   19 &   12 \\
311 & -30935 & -41143 &   19 &   12 \\
311 & -30943 & -46259 &   19 &   12 \\
311 & -30967 & -30967 &   19 &   12 \\
311 & -31015 & -31015 &   19 &   12 \\
311 & -31223 & -31223 &   19 &   12 \\
311 & -31583 & -42007 &   19 &   12 \\
311 & -31655 & -31655 &   19 &   12 \\
311 & -32479 & -48563 &   19 &   12 \\
311 & -32503 & -32503 &   19 &   12 \\
311 & -32563 & -36439 &   19 &   12 \\
311 & -32695 & -32695 &   19 &   12 \\
311 & -32735 & -32735 &   19 &   12 \\
311 & -33203 & -44167 &   19 &   12 \\
311 & -34231 & -34231 &   19 &   12 \\
311 & -34391 & -34391 &   19 &   12 \\
311 & -34535 & -45943 &   19 &   12 \\
311 & -34903 & -34903 &   19 &   12 \\
311 & -35039 & -46615 &   19 &   12 \\
311 & -35047 & -52415 &   19 &   12 \\
311 & -35111 & -35111 &   19 &   12 \\
311 & -35551 & -35551 &   19 &   12 \\
311 & -35831 & -42547 &   19 &   12 \\
311 & -36127 & -54035 &   19 &   12 \\
311 & -36191 & -36191 &   19 &   12 \\
\end{tabular}
\newpage
\noindent\begin{tabular}{rrrrr}
 $k$ & $x_0$ & $x_{lc}$ & $a$ & $b$ \\
\hline
311 & -36319 & -36319 &   19 &   12 \\
311 & -36659 & -48775 &   19 &   12 \\
311 & -37111 & -37111 &   19 &   12 \\
311 & -37235 & -49543 &   19 &   12 \\
311 & -37543 & -37543 &   19 &   12 \\
311 & -37567 & -47191 &   19 &   12 \\
311 & -37747 & -42271 &   19 &   12 \\
311 & -37879 & -37879 &   19 &   12 \\
311 & -37939 & -37939 &   19 &   12 \\
311 & -38119 & -38119 &   19 &   12 \\
311 & -38323 & -38323 &   19 &   12 \\
311 & -38623 & -38623 &   19 &   12 \\
311 & -38999 & -38999 &   19 &   12 \\
311 & -39143 & -39143 &   19 &   12 \\
311 & -39283 & -43999 &   19 &   12 \\
311 & -39359 & -49459 &   19 &   12 \\
311 & -39647 & -39647 &   19 &   12 \\
311 & -39719 & -47155 &   19 &   12 \\
311 & -39775 & -39775 &   19 &   12 \\
311 & -39847 & -39847 &   19 &   12 \\
311 & -40223 & -40223 &   19 &   12 \\
311 & -40279 & -40279 &   19 &   12 \\
311 & -40735 & -40735 &   19 &   12 \\
311 & -40927 & -40927 &   19 &   12 \\
311 & -41303 & -43711 &   19 &   12 \\
311 & -41395 & -41395 &   19 &   12 \\
311 & -41767 & -41767 &   19 &   12 \\
311 & -41843 & -55687 &   19 &   12 \\
311 & -41959 & -59251 &   19 &   12 \\
311 & -42175 & -53023 &   19 &   12 \\
311 & -42355 & -42355 &   19 &   12 \\
311 & -42727 & -42727 &   19 &   12 \\
311 & -42931 & -48103 &   19 &   12 \\
311 & -43283 & -54367 &   19 &   12 \\
311 & -43795 & -49075 &   19 &   12 \\
311 & -44023 & -44023 &   19 &   12 \\
311 & -44255 & -46867 &   19 &   12 \\
311 & -44383 & -44383 &   19 &   12 \\
311 & -44467 & -44467 &   19 &   12 \\
311 & -45427 & -45427 &   19 &   12 \\
311 & -46111 & -46111 &   19 &   12 \\
311 & -46271 & -46271 &   19 &   12 \\
311 & -46451 & -46451 &   19 &   12 \\
\end{tabular}
\hspace{5mm}
\noindent\begin{tabular}{rrrrr}
 $k$ & $x_0$ & $x_{lc}$ & $a$ & $b$ \\
\hline
311 & -47135 & -49855 &   19 &   12 \\
311 & -47551 & -47551 &   19 &   12 \\
311 & -47839 & -47839 &   19 &   12 \\
311 & -47891 & -47891 &   19 &   12 \\
311 & -47911 & -47911 &   19 &   12 \\
311 & -50707 & -56851 &   19 &   12 \\
311 & -51475 & -51475 &   19 &   12 \\
311 & -51571 & -61555 &   19 &   12 \\
311 & -52447 & -52447 &   19 &   12 \\
311 & -53363 & -53363 &   19 &   12 \\
311 & -54643 & -54643 &   19 &   12 \\ \hline
313 &     1 &    35 &  200 &  115 \\
313 &     5 &     5 &   12 &    3 \\
313 &     7 &    47 &   92 &   49 \\
313 &    11 &    11 &   12 &    3 \\
313 & -6329 & -6329 &   68 &   43 \\ \hline
317 &     1 &     1 &   68 &   32 \\
317 &    31 &   539 &   39 &   23 \\ \hline
319 &     1 &    37 &  102 &   54 \\ \hline
323 &     1 &     5 &   58 &   28 \\ \hline
325 &     1 &   127 &   24 &   12 \\
325 &    67 &    83 &   24 &   12 \\ \hline
329 &     1 &     1 &   29 &   10 \\
329 &     5 &    11 &   50 &   22 \\
329 &   179 &   407 &   21 &   12 \\
329 &   295 &   295 &   21 &   12 \\ \hline
331 &     1 &    13 &   52 &   22 \\
331 &     5 &    31 &   63 &   33 \\ \hline
335 &     1 &     1 &   42 &   18 \\
335 &    11 &    17 &   60 &   32 \\
335 &    31 &    73 &   51 &   25 \\ \hline
337 &     1 &    25 &   32 &   16 \\
337 &     5 &    23 &   32 &   16 \\
337 &    13 &    47 &   16 &    8 \\
337 &    43 &   227 &   16 &    8 \\
337 &    59 &    59 &   32 &   16 \\
337 &    91 &    91 &   16 &    8 \\
337 &   107 &   107 &   16 &    8 \\
337 &   209 &   209 &   16 &    8 \\ \hline
341 &     1 &  1003 &   50 &   30 \\
341 &   371 &   647 &   50 &   30 \\ \hline
343 &     1 &   221 &   38 &   22 \\
343 &     5 &    13 &   41 &   16 \\
\end{tabular}
\newpage
\noindent\begin{tabular}{rrrrr}
 $k$ & $x_0$ & $x_{lc}$ & $a$ & $b$ \\
\hline
343 & 177337 & 177337 &  149 &   94 \\ \hline
347 &     1 &     7 &   62 &   30 \\
347 &     5 &     5 &  160 &   83 \\ \hline
349 &     1 &    11 &   58 &   29 \\
349 &     5 &    31 &  116 &   58 \\ \hline
353 &     1 &     1 &   71 &   28 \\
353 &     5 &   191 &   28 &   16 \\
353 &    25 &   265 &   28 &   16 \\
353 &    71 &   167 &   56 &   32 \\
353 &    95 &   319 &   28 &   16 \\ \hline
355 &     1 &     1 &   13 &    3 \\
355 &     7 &     7 &   13 &    3 \\
355 &    49 &    49 &   40 &   20 \\
355 &  4573 & 14429 &   27 &   17 \\
355 &  6037 & 18413 &   27 &   17 \\
355 &  7021 & 13781 &   27 &   17 \\
355 &  7409 & 14077 &   27 &   17 \\
355 &  8669 & 13181 &   27 &   17 \\
355 &  9841 & 11293 &   27 &   17 \\
355 & 11221 & 17009 &   27 &   17 \\
355 & 11377 & 11377 &   81 &   51 \\
355 & 12181 & 23821 &   27 &   17 \\
355 & 13933 & 21077 &   27 &   17 \\
355 & 15581 & 15581 &   27 &   17 \\
355 & 15997 & 15997 &   27 &   17 \\
355 & 18461 & 18461 &   27 &   17 \\
355 & 20693 & 27709 &   27 &   17 \\
355 & 22909 & 22909 &   27 &   17 \\
355 & 26221 & 26221 &   27 &   17 \\
355 & 30449 & 30449 &   27 &   17 \\ \hline
359 &     1 &   169 &   15 &    8 \\
359 &     7 &   113 &   15 &    8 \\
359 &    13 &   269 &   15 &    8 \\
359 &    17 &   121 &   15 &    8 \\
359 &    41 &   229 &   15 &    8 \\
359 &    89 &   145 &   15 &    8 \\
359 &   109 &   109 &   30 &   16 \\
359 &   157 &   157 &   15 &    8 \\
359 &   191 &   191 &   15 &    8 \\
359 &   209 &   209 &   15 &    8 \\
359 &   263 &   263 &   15 &    8 \\ \hline
361 &     1 &   113 &   11 &    5 \\
361 &     5 &     5 &   55 &   25 \\
\end{tabular}
\hspace{5mm}
\noindent\begin{tabular}{rrrrr}
 $k$ & $x_0$ & $x_{lc}$ & $a$ & $b$ \\
\hline
361 &     7 &    53 &   22 &   10 \\
361 &    11 &    59 &   22 &   10 \\
361 &    13 &   137 &   11 &    5 \\
361 &    17 &    17 &   22 &   10 \\
361 &    35 &    35 &   22 &   10 \\
361 &    67 &    79 &   11 &    5 \\
361 &    97 &    97 &   11 &    5 \\
361 &   131 &   131 &   11 &    5 \\ \hline
365 &     1 &   127 &   96 &   56 \\
365 &     7 &    11 &   12 &    4 \\
365 &    13 &    13 &   12 &    4 \\
365 &    19 &    41 &   12 &    4 \\ \hline
367 &     1 &     7 &  135 &   62 \\ \hline
371 &     1 &    25 &  222 &  114 \\
371 &  -499 &  -563 &   14 &   10 \\ \hline
373 &     1 &    31 &   86 &   41 \\
373 &     5 &    23 &   86 &   41 \\ \hline
377 &     1 &    19 &  124 &   64 \\
377 &    17 &    53 &   44 &   20 \\ \hline
379 &     1 &    61 &   28 &   14 \\
379 &    35 &    35 &   28 &   14 \\
379 &    41 &    41 &   28 &   14 \\ \hline
383 &     1 &    67 &   81 &   43 \\
383 &     5 &    41 &   44 &   21 \\ \hline
385 &     1 &   107 &   84 &   48 \\
385 &    13 &    17 &   18 &    6 \\
385 &    23 &    23 &   40 &   20 \\
385 & -1409 & -3457 &   22 &   14 \\
385 & -4241 & -4241 &   22 &   14 \\
385 & -4321 & -5221 &   22 &   14 \\
385 & -5473 & -5473 &   22 &   14 \\ \hline
389 &     1 &    59 &   87 &   49 \\
389 &    17 &    17 &   31 &   13 \\
389 &    23 &    37 &   87 &   49 \\ \hline
391 &     1 &    77 &   63 &   34 \\
391 &     5 &    37 &   20 &    8 \\
391 &    13 &    13 &   20 &    8 \\
391 &    29 &    29 &   20 &    8 \\ \hline
395 &     1 &    83 &   22 &   10 \\
395 &     7 &    47 &   22 &   10 \\
395 &    17 &    17 &   22 &   10 \\
395 &    19 &   139 &   22 &   10 \\
395 &    23 &    29 &   22 &   10 \\ \hline
\end{tabular}
\newpage
\noindent\begin{tabular}{rrrrr}
 $k$ & $x_0$ & $x_{lc}$ & $a$ & $b$ \\
\hline
397 &     1 &    11 &  104 &   54 \\ \hline
401 &     1 &     5 &  188 &   88 \\ \hline
403 &     1 &     5 &  156 &   84 \\ \hline
407 &     1 &     1 &   46 &   17 \\
407 &     5 &    31 &   90 &   45 \\
407 &    13 &    71 &   74 &   43 \\
407 &    23 &   119 &   30 &   15 \\
407 & -3095 & -3095 &   44 &   28 \\
407 & -3923 & -6035 &   44 &   28 \\ \hline
409 &     1 &    23 &   49 &   23 \\
409 &     7 &    19 &   49 &   23 \\
409 &   353 &   371 &   38 &   22 \\ \hline
413 &     1 &     5 &   86 &   40 \\
413 &   515 &   979 &   30 &   18 \\
413 &   659 &   911 &   30 &   18 \\ \hline
415 &     1 &    23 &  140 &   68 \\ \hline
419 &     1 &    41 &  186 &   98 \\
419 &   557 &  1045 &   28 &   17 \\ \hline
421 &     1 &    37 &   52 &   23 \\
421 &     5 &    67 &   88 &   47 \\ \hline
425 &     1 &   137 &   32 &   16 \\
425 &    19 &    73 &   32 &   16 \\
425 &    23 &    23 &   32 &   16 \\
425 &    59 &    83 &   32 &   16 \\ \hline
427 &     1 &    67 &  120 &   60 \\ \hline
431 &     1 &   593 &   32 &   19 \\
431 &     5 &   157 &    9 &    4 \\
431 &     7 &     7 &   27 &   12 \\
431 &    11 &   151 &    9 &    4 \\
431 &    13 &    37 &   18 &    8 \\
431 &    17 &    73 &    9 &    4 \\
431 &    19 &    19 &   59 &   31 \\
431 &    23 &   125 &    9 &    4 \\
431 &    25 &    25 &   18 &    8 \\
431 &    31 &   103 &    9 &    4 \\
431 &    35 &    35 &   36 &   16 \\
431 &    47 &   143 &    9 &    4 \\
431 &    53 &    53 &   18 &    8 \\
431 &    65 &    65 &    9 &    4 \\
431 &    77 &    89 &    9 &    4 \\
431 &    83 &    85 &    9 &    4 \\
431 &    95 &   121 &    9 &    4 \\
431 &   101 &   101 &    9 &    4 \\
\end{tabular}
\hspace{5mm}
\noindent\begin{tabular}{rrrrr}
 $k$ & $x_0$ & $x_{lc}$ & $a$ & $b$ \\
\hline
431 &   119 &   119 &    9 &    4 \\
431 &   133 &   133 &    9 &    4 \\
431 &   175 &   175 &    9 &    4 \\
431 &   353 &   745 &   32 &   19 \\
431 &   737 &   737 &   32 &   19 \\
431 & -1103 & -1775 &   23 &   15 \\ \hline
433 &     1 &    13 &  208 &  114 \\
433 &     5 &    31 &  128 &   66 \\
433 &    61 &    61 &   24 &    9 \\
433 &   491 &   799 &   40 &   24 \\ \hline
437 &     1 &     1 &  121 &   55 \\
437 &     5 &     5 &   30 &   12 \\
437 &   331 &   587 &   31 &   19 \\ \hline
439 &     1 &   149 &   13 &    6 \\
439 &     5 &     5 &   78 &   36 \\
439 &     7 &   151 &   13 &    6 \\
439 &    11 &    95 &   26 &   12 \\
439 &    25 &    85 &   13 &    6 \\
439 &    41 &    41 &   13 &    6 \\
439 &    47 &   125 &   13 &    6 \\
439 &    61 &    61 &   26 &   12 \\
439 &    65 &    65 &   13 &    6 \\
439 &  4513 &  4513 &  200 &  126 \\ \hline
443 &     1 &    91 &   16 &    7 \\
443 &     7 &    59 &   16 &    7 \\
443 &    13 &    13 &   32 &   14 \\
443 &    19 &    19 &   32 &   14 \\
443 &    41 &    41 &   16 &    7 \\ \hline
445 &     1 &     1 &   28 &    8 \\
445 &    11 &    29 &  108 &   56 \\ \hline
449 &     1 &    19 &   53 &   22 \\ \hline
451 &     1 &    43 &   80 &   40 \\
451 & -1139 & -1483 &   30 &   20 \\ \hline
455 &     1 &    31 &   48 &   24 \\
455 &    11 &    43 &   48 &   24 \\
455 &    17 &    17 &   84 &   48 \\
455 &    37 &   103 &   24 &   12 \\ \hline
457 &     1 &    49 &   69 &   33 \\
457 &     5 &   287 &   37 &   21 \\
457 &    13 &    13 &   32 &   12 \\
457 &   103 &   103 &   37 &   21 \\ \hline
461 &     1 &     7 &   47 &   17 \\ \hline
463 &     1 &    29 &  138 &   72 \\ \hline
\end{tabular}
\newpage
\noindent\begin{tabular}{rrrrr}
 $k$ & $x_0$ & $x_{lc}$ & $a$ & $b$ \\
\hline
467 &     1 &   107 &   50 &   26 \\
467 &     5 &    71 &   50 &   26 \\
467 &     7 &   125 &   50 &   26 \\
467 & -436099 & -436099 &   84 &   53 \\
467 & -460571 & -613939 &   84 &   53 \\
467 & -479039 & -638563 &   84 &   53 \\
467 & -493555 & -740099 &   84 &   53 \\
467 & -507839 & -507839 &   84 &   53 \\
467 & -580595 & -580595 &   84 &   53 \\
467 & -612383 & -979835 &   84 &   53 \\
467 & -663811 & -995483 &   84 &   53 \\
467 & -697411 & -1045883 &   84 &   53 \\
467 & -707071 & -1006079 &   84 &   53 \\
467 & -714523 & -772691 &   84 &   53 \\
467 & -715507 & -966299 &   84 &   53 \\
467 & -721787 & -760699 &   84 &   53 \\
467 & -726847 & -726847 &   84 &   53 \\
467 & -1133123 & -1133123 &   84 &   53 \\
467 & -1420291 & -1457663 &   84 &   53 \\ \hline
469 &     1 &    31 &   84 &   36 \\ \hline
473 &     1 &     1 &  182 &   84 \\
473 &  -569 &  -665 &    8 &    6 \\
473 &  -697 &  -697 &    8 &    6 \\
473 &  -745 &  -745 &    8 &    6 \\ \hline
475 &     1 &    11 &   38 &   14 \\
475 &     7 &    29 &   72 &   36 \\
475 & -1243 & -2011 &   17 &   11 \\
475 & -1739 & -2459 &   17 &   11 \\
475 & -1787 & -2443 &   17 &   11 \\
475 & -2251 & -3431 &   17 &   11 \\
475 & -2267 & -2267 &   17 &   11 \\
475 & -2531 & -2531 &   17 &   11 \\
475 & -2663 & -2699 &   17 &   11 \\
475 & -3107 & -3107 &   17 &   11 \\ \hline
479 &     1 &     7 &  105 &   52 \\
479 &     5 &    49 &  174 &   93 \\
479 &    13 &    37 &   47 &   21 \\
479 &    29 &    41 &   58 &   31 \\ \hline
481 &     1 &   119 &   24 &   12 \\
481 &     5 &   155 &   12 &    6 \\
481 &     7 &   113 &   24 &   12 \\
481 &    11 &   175 &   12 &    6 \\
481 &    23 &    23 &   36 &   18 \\
\end{tabular}
\hspace{5mm}
\noindent\begin{tabular}{rrrrr}
 $k$ & $x_0$ & $x_{lc}$ & $a$ & $b$ \\
\hline
481 &    25 &   139 &   24 &   12 \\
481 &    41 &   215 &   12 &    6 \\
481 &    95 &    95 &   12 &    6 \\
481 &   127 &   127 &   12 &    6 \\
481 &   217 &   217 &   12 &    6 \\
481 &   233 &   233 &   12 &    6 \\
481 &   281 &   281 &   12 &    6 \\
481 &   353 &   353 &   12 &    6 \\ \hline
485 &     1 &   283 &   43 &   25 \\
485 &     7 &     7 &  112 &   64 \\
485 &    11 &   259 &   26 &   14 \\
485 &    13 &    83 &   35 &   17 \\
485 &    17 &   251 &   26 &   14 \\
485 &    19 &    19 &    9 &    3 \\
485 &    23 &    23 &    9 &    3 \\
485 &    29 &    29 &    9 &    3 \\
485 &    31 &    31 &    9 &    3 \\
485 &    37 &    37 &    9 &    3 \\
485 &    47 &    47 &    9 &    3 \\
485 &    49 &    49 &    9 &    3 \\
485 &    53 &    53 &    9 &    3 \\
485 &   101 &   101 &   44 &   20 \\
485 &   103 &   299 &   43 &   25 \\
485 &  -869 & -2021 &   17 &   11 \\
485 & -1733 & -2149 &   17 &   11 \\
485 & -2293 & -3197 &   17 &   11 \\
485 & -2341 & -2341 &   17 &   11 \\
485 & -2597 & -2597 &   17 &   11 \\
485 & -3193 & -3289 &   17 &   11 \\
485 & -3673 & -3673 &   17 &   11 \\
485 & -4249 & -4249 &   17 &   11 \\ \hline
487 &     1 &   103 &   88 &   46 \\
487 &    11 &    43 &  131 &   74 \\ \hline
491 &     1 &    49 &  110 &   57 \\
491 &    11 &    11 &   80 &   37 \\
491 & -1003 & -1003 &   30 &   20 \\ \hline
493 &     1 &     1 &  200 &   96 \\
493 &  -749 & -1261 &   12 &    8 \\
493 &  -973 & -1453 &   12 &    8 \\
493 & -1229 & -1597 &   12 &    8 \\
493 & -1325 & -1325 &   12 &    8 \\
493 & -1405 & -1709 &   12 &    8 \\
493 & -1469 & -1469 &   12 &    8 \\
\end{tabular}
\newpage
\noindent\begin{tabular}{rrrrr}
 $k$ & $x_0$ & $x_{lc}$ & $a$ & $b$ \\
\hline
493 & -1613 & -1613 &   12 &    8 \\
493 & -1637 & -2141 &   12 &    8 \\
493 & -1829 & -1829 &   12 &    8 \\ \hline
497 &     1 &   103 &   40 &   20 \\
497 &    11 &   265 &   20 &   10 \\
497 &    13 &    13 &  140 &   70 \\
497 &    19 &   145 &   20 &   10 \\
497 &    55 &   113 &   20 &   10 \\
497 &    79 &    79 &   20 &   10 \\
497 &    95 &    95 &   20 &   10 \\ \hline
499 &     1 &     1 &   47 &   19 \\
499 &     5 &   139 &  191 &  109 \\
499 &    11 &    23 &   47 &   19 \\
499 &   613 &  2381 &   48 &   30 \\
499 &   901 &  1549 &   80 &   50 \\
499 &  1057 &  5501 &   32 &   20 \\
499 &  1261 &  4781 &   32 &   20 \\
499 &  2029 &  7085 &   16 &   10 \\
499 &  3053 &  6605 &   16 &   10 \\
499 &  3269 &  4525 &   32 &   20 \\
499 &  3469 &  5453 &   16 &   10 \\
499 &  3649 &  5281 &   32 &   20 \\
499 &  3853 &  6029 &   16 &   10 \\
499 &  3917 &  4513 &   32 &   20 \\
499 &  4061 &  5581 &   16 &   10 \\
499 &  4189 &  4189 &   32 &   20 \\
499 &  4205 &  5773 &   16 &   10 \\
499 &  4429 &  6893 &   16 &   10 \\
499 &  4549 &  7073 &   16 &   10 \\
499 &  4573 &  7109 &   16 &   10 \\
499 &  4997 & 10465 &   16 &   10 \\
499 &  5005 &  5005 &   16 &   10 \\
499 &  5101 &  7901 &   16 &   10 \\
499 &  5213 &  7117 &   16 &   10 \\
499 &  5317 &  8225 &   16 &   10 \\
499 &  5357 &  5357 &   16 &   10 \\
499 &  5509 &  8513 &   16 &   10 \\
499 &  5533 &  8549 &   16 &   10 \\
499 &  5725 &  8837 &   16 &   10 \\
499 &  5965 &  5965 &   16 &   10 \\
499 &  6085 &  6085 &   16 &   10 \\
499 &  6221 &  6221 &   16 &   10 \\
499 &  6365 &  6365 &   16 &   10 \\
\end{tabular}
\hspace{5mm}
\noindent\begin{tabular}{rrrrr}
 $k$ & $x_0$ & $x_{lc}$ & $a$ & $b$ \\
\hline
499 &  6509 &  7585 &   16 &   10 \\
499 &  6689 &  9085 &   16 &   10 \\
499 &  6725 &  7841 &   16 &   10 \\
499 &  6853 &  6853 &   16 &   10 \\
499 &  7013 &  9445 &   16 &   10 \\
499 &  7045 & 10817 &   16 &   10 \\
499 &  7069 &  7069 &   16 &   10 \\
499 &  7201 &  7201 &   16 &   10 \\
499 &  7261 &  7261 &   16 &   10 \\
499 &  7453 &  7517 &   16 &   10 \\
499 &  7489 &  8737 &   16 &   10 \\
499 &  7525 &  7525 &   16 &   10 \\
499 &  7661 &  8581 &   16 &   10 \\
499 &  7777 &  9061 &   16 &   10 \\
499 &  7789 &  7789 &   16 &   10 \\
499 &  8545 &  8545 &   16 &   10 \\
499 &  8797 &  8797 &   16 &   10 \\
499 & 10433 & 10433 &   16 &   10 \\
499 & 10849 & 10849 &   16 &   10 \\ \hline
503 &     1 &    11 &    9 &    2 \\
503 &     5 &     5 &    9 &    2 \\
503 &     7 &     7 &    9 &    2 \\
503 &    13 &    23 &   86 &   47 \\
503 &    17 &    29 &  104 &   51 \\
503 &    19 &    19 &    9 &    2 \\
503 &   173 &   509 &   77 &   45 \\
503 & -15391 & -24607 &   68 &   43 \\ \hline
505 &     1 &   113 &   14 &    6 \\
505 &     7 &    79 &   14 &    6 \\
505 &    11 &    91 &   14 &    6 \\
505 &    13 &    59 &   14 &    6 \\
505 &    31 &    49 &   14 &    6 \\ \hline
509 &     1 &     1 &    9 &    1 \\
509 &     5 &     5 &  262 &  142 \\ \hline
511 &     1 &    19 &   72 &   36 \\
511 &     5 &    11 &   54 &   24 \\
511 &    13 &   103 &   27 &   12 \\
511 &    23 &    23 &   27 &   12 \\
511 &    25 &   293 &   45 &   24 \\ \hline
515 &     1 &     7 &   66 &   30 \\
515 &    11 &    43 &   66 &   30 \\ \hline
517 &     1 &   379 &   14 &    8 \\
517 &     5 &   125 &   70 &   40 \\
\end{tabular}
\newpage
\noindent\begin{tabular}{rrrrr}
 $k$ & $x_0$ & $x_{lc}$ & $a$ & $b$ \\
\hline
517 &     7 &   251 &   42 &   24 \\
517 &    17 &   173 &   42 &   24 \\
517 &    23 &   443 &   14 &    8 \\
517 &   227 &   787 &   14 &    8 \\
517 &   263 &   523 &   14 &    8 \\
517 &   299 &   571 &   14 &    8 \\
517 &   397 &   427 &   14 &    8 \\
517 &   487 &   871 &   14 &    8 \\
517 &   515 &   515 &   14 &    8 \\
517 &   643 &   643 &   14 &    8 \\ \hline
521 &     1 &    55 &   95 &   50 \\ \hline
523 &     1 &    89 &   45 &   23 \\
523 &     5 &   115 &   45 &   23 \\
523 &    13 &   209 &   45 &   23 \\ \hline
527 &     1 &    35 &  101 &   54 \\
527 &     5 &    71 &   65 &   30 \\
527 &    13 &    13 &   47 &   18 \\
527 & -1199 & -1199 &   18 &   12 \\
527 & -1279 & -1655 &   18 &   12 \\ \hline
529 &     1 &    29 &   42 &   19 \\
529 &     5 &   241 &   53 &   30 \\
529 &    13 &    49 &   42 &   19 \\ \hline
533 &     1 &    59 &   20 &    8 \\ \hline
535 &     1 &     7 &  208 &  112 \\
535 &    19 &    47 &   78 &   42 \\
535 &    79 &   193 &   26 &   14 \\ \hline
539 &     1 &    19 &   18 &    6 \\
539 &    31 &    31 &   18 &    6 \\ \hline
541 &     1 &    25 &  180 &   90 \\ \hline
545 &     1 &    31 &   28 &   12 \\
545 &    49 &    49 &   28 &   12 \\ \hline
547 &     1 &   283 &   39 &   23 \\
547 &     5 &    13 &  156 &   78 \\
547 &   257 &   305 &   39 &   23 \\ \hline
551 &     1 &    37 &   44 &   20 \\
551 &    11 &   101 &   22 &   10 \\
551 &    47 &    59 &   44 &   20 \\
551 & 78257 & 199345 &   65 &   41 \\
551 & 105841 & 141305 &   65 &   41 \\
551 & 120665 & 120665 &   65 &   41 \\
551 & 132977 & 132977 &   65 &   41 \\
551 & 141593 & 212665 &   65 &   41 \\
551 & 146489 & 146489 &   65 &   41 \\
\end{tabular}
\hspace{5mm}
\noindent\begin{tabular}{rrrrr}
 $k$ & $x_0$ & $x_{lc}$ & $a$ & $b$ \\
\hline
551 & 151345 & 151345 &   65 &   41 \\ \hline
553 &     1 &     1 &   21 &    6 \\
553 &     5 &    13 &   66 &   30 \\
553 &    11 &   239 &   45 &   24 \\ \hline
557 &     1 &    23 &   23 &    9 \\
557 &     5 &     5 &   46 &   18 \\
557 &     7 &    47 &   23 &    9 \\
557 &    17 &    17 &   23 &    9 \\ \hline
559 &     1 &     1 &   60 &   24 \\ \hline
563 &     1 &    19 &  426 &  215 \\ \hline
565 &     1 &    13 &  244 &  128 \\
565 &     7 &    19 &   68 &   32 \\ \hline
569 &     1 &    11 &  136 &   64 \\ \hline
571 &     1 &    19 &  112 &   50 \\
571 &    17 &    61 &   67 &   35 \\ \hline
575 &     1 &    31 &   88 &   44 \\
575 &     7 &    19 &   47 &   21 \\ \hline
577 &     1 &    91 &   15 &    7 \\
577 &     5 &   217 &   15 &    7 \\
577 &    11 &    17 &   45 &   21 \\
577 &    13 &    79 &   15 &    7 \\
577 &    31 &   107 &   15 &    7 \\
577 &    41 &    41 &   30 &   14 \\
577 &    47 &    47 &   15 &    7 \\
577 &    61 &    61 &   30 &   14 \\
577 &   121 &   151 &   15 &    7 \\
577 &   155 &   155 &   15 &    7 \\ \hline
581 &     1 &    53 &  118 &   62 \\ \hline
583 &     1 &     7 &   76 &   32 \\ \hline
587 &     1 &   181 &   30 &   16 \\
587 &    35 &   173 &   60 &   32 \\
587 &   203 &   299 &   30 &   16 \\ \hline
589 &     1 &    29 &  126 &   54 \\
589 &     5 &    91 &   78 &   42 \\ \hline
593 &     1 &     1 &  108 &   48 \\
593 &     5 &    29 &   54 &   24 \\
593 &    23 &    23 &   18 &    8 \\
593 &    71 &    77 &   18 &    8 \\ \hline
595 &     1 &   137 &   32 &   16 \\
595 &    11 &    13 &   64 &   32 \\
595 &    29 &    29 &   32 &   16 \\ \hline
599 &     1 &   205 &   36 &   19 \\
599 &     5 &    95 &   36 &   19 \\
\end{tabular}
\newpage
\noindent\begin{tabular}{rrrrr}
 $k$ & $x_0$ & $x_{lc}$ & $a$ & $b$ \\
\hline
599 &     7 &   175 &   36 &   19 \\ \hline
601 &     1 &   211 &   30 &   15 \\
601 &     5 &     5 &   16 &    3 \\
601 &     7 &     7 &   23 &    9 \\
601 &    11 &   103 &   30 &   15 \\
601 &    23 &  1303 &   44 &   27 \\
601 &    25 &   313 &   81 &   48 \\
601 &    29 &   161 &   74 &   42 \\
601 &    31 &    31 &   23 &    9 \\
601 &    49 &   131 &   53 &   24 \\
601 &    65 &    65 &   30 &   15 \\
601 &   115 &   115 &   30 &   15 \\
601 &   119 &   119 &   74 &   42 \\
601 &  -665 &  -665 &    7 &    6 \\ \hline
605 &     1 &    43 &  198 &  106 \\
605 &    13 &    13 &  110 &   50 \\
605 & -2141 & -5213 &   44 &   28 \\
605 & -3121 & -4397 &   66 &   42 \\
605 & -3901 & -5933 &   22 &   14 \\
605 & -5597 & -7261 &   22 &   14 \\
605 & -5677 & -5677 &   44 &   28 \\
605 & -6109 & -6109 &   22 &   14 \\
605 & -6893 & -8989 &   22 &   14 \\
605 & -6973 & -6973 &   22 &   14 \\
605 & -7217 & -9421 &   22 &   14 \\
605 & -7597 & -11093 &   22 &   14 \\
605 & -7981 & -7981 &   22 &   14 \\
605 & -8341 & -8341 &   22 &   14 \\
605 & -8621 & -8621 &   22 &   14 \\
605 & -9773 & -9773 &   22 &   14 \\
605 & -10033 & -10033 &   22 &   14 \\
605 & -10837 & -14389 &   22 &   14 \\ \hline
607 &     1 &    23 &   27 &   12 \\
607 &    13 &    13 &   27 &   12 \\
607 &    25 &    25 &   27 &   12 \\
607 &    47 &    85 &   27 &   12 \\
607 &    65 &    65 &   27 &   12 \\ \hline
611 &     1 &   331 &   56 &   32 \\
611 &     5 &   349 &   28 &   16 \\
611 &    23 &   187 &   84 &   48 \\
611 &   173 &   629 &   28 &   16 \\
611 &   515 &   515 &   28 &   16 \\ \hline
613 &     1 &     5 &   30 &   11 \\
\end{tabular}
\hspace{5mm}
\noindent\begin{tabular}{rrrrr}
 $k$ & $x_0$ & $x_{lc}$ & $a$ & $b$ \\
\hline
613 &    11 &    11 &   30 &   11 \\
613 &  4555 &  7139 &   78 &   49 \\ \hline
617 &     1 &    19 &  147 &   84 \\
617 &     7 &   547 &   49 &   28 \\
617 &    53 &   383 &   49 &   28 \\ \hline
619 &     1 &    67 &  177 &   93 \\ \hline
623 &     1 &   145 &   74 &   40 \\
623 &     5 &    43 &   79 &   32 \\ \hline
625 &     1 &    29 &  101 &   43 \\
625 &    59 &   569 &   33 &   19 \\ \hline
629 &     1 &     1 &   28 &    8 \\
629 &     5 &    59 &  164 &   88 \\
629 &   535 &   787 &   68 &   40 \\ \hline
631 &     1 &   353 &   97 &   56 \\ \hline
635 &     1 &    19 &  148 &   72 \\
635 &   893 &   893 &   30 &   18 \\ \hline
637 &     1 &     1 &  156 &   66 \\ \hline
641 &     1 &     7 &  214 &  100 \\
641 & -1025 & -1153 &   30 &   20 \\
641 & -1537 & -1537 &   15 &   10 \\
641 & -1633 & -2129 &   15 &   10 \\
641 & -1729 & -1729 &   15 &   10 \\
641 & -1817 & -1817 &   15 &   10 \\
641 & -1945 & -2521 &   15 &   10 \\
641 & -2017 & -2017 &   15 &   10 \\ \hline
643 &     1 &     1 &   19 &    5 \\
643 &     7 &     7 &   19 &    5 \\
643 &    11 &    37 &   56 &   26 \\
643 &    13 &    19 &   93 &   47 \\
643 &    25 &    77 &  148 &   84 \\
643 &   379 &   379 &   37 &   21 \\ \hline
647 &     1 &     7 &  174 &   91 \\
647 &     5 &    65 &   59 &   29 \\
647 &   113 &   493 &   56 &   33 \\ \hline
649 &     1 &    25 &  134 &   63 \\
649 &     7 &    97 &   78 &   41 \\
649 &   707 &  1063 &   50 &   30 \\ \hline
653 &     1 &     5 &  223 &  105 \\ \hline
655 &     1 &     7 &  194 &  102 \\ \hline
659 &     1 &   293 &   44 &   25 \\
659 &    19 &    71 &   88 &   50 \\ \hline
661 &     1 &   455 &   60 &   36 \\
661 &     5 &    17 &  180 &   86 \\ \hline
\end{tabular}
\newpage
\noindent\begin{tabular}{rrrrr}
 $k$ & $x_0$ & $x_{lc}$ & $a$ & $b$ \\
\hline
665 &     1 &   139 &   78 &   42 \\
665 &    17 &    23 &   66 &   30 \\ \hline
667 &     1 &   443 &   17 &    9 \\
667 &     5 &   221 &   17 &    9 \\
667 &     7 &   293 &   17 &    9 \\
667 &    13 &   211 &   17 &    9 \\
667 &    19 &   181 &   17 &    9 \\
667 &    67 &   329 &   17 &    9 \\
667 &   157 &   569 &   17 &    9 \\
667 &   275 &   275 &   17 &    9 \\
667 &   281 &   281 &   17 &    9 \\
667 &   635 &   635 &   17 &    9 \\ \hline
671 &     1 &   151 &  120 &   70 \\
671 &    13 &   295 &   30 &   15 \\
671 &    25 &   169 &   30 &   15 \\
671 &    37 &    89 &   30 &   15 \\
671 &    71 &   167 &   30 &   15 \\
671 &   161 &   161 &   30 &   15 \\ \hline
673 &     1 &     1 &   22 &    7 \\
673 &     5 &    31 &  146 &   77 \\
673 &     7 &    13 &   22 &    7 \\ \hline
677 &     1 &   301 &   31 &   17 \\
677 &     5 &   119 &   31 &   17 \\ \hline
679 &     1 &   103 &   96 &   48 \\
679 &     5 &    59 &   96 &   48 \\
679 &    19 &    23 &   18 &    6 \\
679 &  1009 &  1009 &   60 &   36 \\ \hline
683 &     1 &    35 &  176 &   93 \\ \hline
685 &     1 &    13 &  128 &   56 \\
685 &    31 &    31 &   72 &   40 \\ \hline
689 &     1 &    79 &  148 &   76 \\ \hline
691 &     1 &    31 &  111 &   51 \\ \hline
695 &     1 &    31 &  125 &   67 \\
695 &     7 &    13 &   81 &   39 \\
695 & -2743 & -6839 &   22 &   14 \\
695 & -5431 & -7799 &   22 &   14 \\
695 & -6439 & -9311 &   22 &   14 \\
695 & -6967 & -10103 &   22 &   14 \\
695 & -7519 & -10931 &   22 &   14 \\
695 & -10547 & -13831 &   22 &   14 \\
695 & -10591 & -10591 &   22 &   14 \\
695 & -11383 & -11383 &   22 &   14 \\
695 & -11923 & -11923 &   22 &   14 \\ \hline
\end{tabular}
\hspace{5mm}
\noindent\begin{tabular}{rrrrr}
 $k$ & $x_0$ & $x_{lc}$ & $a$ & $b$ \\
\hline
697 &     1 &     1 &   85 &   38 \\
697 &     5 &   263 &   25 &   14 \\
697 &     7 &     7 &   20 &    8 \\
697 &    19 &    19 &   40 &   16 \\
697 &    25 &    25 &   20 &    8 \\
697 &    29 &   665 &   25 &   14 \\
697 &    31 &    47 &   20 &    8 \\
697 &    37 &    67 &   20 &    8 \\
697 &    43 &   227 &   75 &   42 \\
697 &    79 &    79 &   45 &   22 \\
697 &   199 &   647 &   25 &   14 \\
697 & -4433 & -6073 &   85 &   54 \\ \hline
701 &     1 &    83 &  111 &   59 \\ \hline
703 &     1 &    71 &   36 &   18 \\
703 &    49 &    97 &   36 &   18 \\ \hline
707 &     1 &    47 &   42 &   18 \\
707 &    11 &    25 &   84 &   36 \\ \hline
709 &     1 &   203 &  188 &  104 \\ \hline
713 &     1 &     1 &   79 &   36 \\
713 &     5 &   401 &   91 &   54 \\ \hline
715 &     1 &     7 &   36 &   12 \\
715 &    17 &    59 &  112 &   64 \\
715 &    53 &    53 &   56 &   32 \\
715 &    73 &   389 &   28 &   16 \\ \hline
719 &     1 &     1 &  151 &   76 \\
719 &     7 &   583 &   37 &   21 \\
719 &    11 &   215 &   74 &   42 \\
719 &    31 &   433 &   37 &   21 \\
719 &    35 &   625 &   37 &   21 \\
719 &   209 &   673 &   37 &   21 \\ \hline
721 &     1 &   145 &   69 &   36 \\
721 &     5 &    73 &  138 &   72 \\
721 &    11 &    11 &   69 &   36 \\ \hline
725 &     1 &    49 &   54 &   22 \\
725 &     7 &    11 &  112 &   56 \\
725 &    23 &    31 &   54 &   22 \\
725 &    43 &   157 &   58 &   34 \\ \hline
727 &     1 &    79 &   44 &   22 \\
727 &     5 &    89 &   44 &   22 \\
727 &    11 &  1577 &   45 &   28 \\
727 &    13 &    49 &   44 &   22 \\
727 &    53 &   125 &   44 &   22 \\ \hline
731 &     1 &    35 &   60 &   24 \\
\end{tabular}%

\section{Limit cycles: k from 1001 to 1999}

\noindent%

\section{\bf Limit cycles: k from 2003 to 2999}

\noindent%

\section{Limit cycles: k from 3001 to 3997}

\noindent%

\section{\bf Limit cycles: k from 4001 to 4999}

\noindent%

\section{Limit cycles: k from 5003 to 5999}

\noindent%

\section{Limit cycles: k from 6001 to 6997}

\noindent%

\section{Limit cycles: k from 7001 to 7999}

\noindent%

\section{Limit cycles: k from 8003 to 8999}

\noindent%

\section{Limit cycles: k from 9001 to 9997}

\noindent%
